\documentclass[reqno,10pt]{amsart}
\usepackage{xcolor}

\usepackage{amssymb}
\usepackage{comment}
\usepackage{bbm}
\usepackage{tikz}
\usepackage[all,cmtip]{xy}
\usepackage{hyperref}
\usepackage{stmaryrd}

\newcounter{mt}

\newtheorem{MainTheorem}[mt]{Theorem}

\newtheorem{Proposition}{Proposition}[section]
\newtheorem{Definition}[Proposition]{Definition}
\newtheorem{Lemma}[Proposition]{Lemma}
\newtheorem{Theorem}[Proposition]{Theorem}

\DeclareMathOperator{\Val}{Val}
\DeclareMathOperator{\Curv}{Curv}

\DeclareMathOperator{\nc}{nc}

\DeclareMathOperator{\tr}{tr}

\newcommand{\p}{\mathbb{P}}
\DeclareMathOperator{\vol}{vol}

\DeclareMathOperator{\Dens}{Dens}

\DeclareMathOperator{\Stab}{Stab}
\DeclareMathOperator{\Sym}{Sym}

\DeclareMathOperator{\WF}{WF}
\DeclareMathOperator{\LC}{LC}

\DeclareMathOperator{\OO}{O}
\DeclareMathOperator{\SO}{SO}

\DeclareMathOperator{\Isom}{Isom}
\newcommand{\R}{\mathbb{R}}
\newcommand{\C}{\mathbb{C}}

\makeatletter
\def\moverlay{\mathpalette\mov@rlay}
\def\mov@rlay#1#2{\leavevmode\vtop{%
		\baselineskip\z@skip \lineskiplimit-\maxdimen
		\ialign{\hfil$\m@th#1##$\hfil\cr#2\crcr}}}
\newcommand{\charfusion}[3][\mathord]{
	#1{\ifx#1\mathop\vphantom{#2}\fi
		\mathpalette\mov@rlay{#2\cr#3}
	}
	\ifx#1\mathop\expandafter\displaylimits\fi}
\makeatother

\newcommand{\largewedge}{\mbox{\Large $\wedge$}}

\newtheorem{note}{Remark} 
\def\note#1{\ifvmode\leavevmode\fi\vadjust{\vbox to0pt{\vss
 \hbox to 0pt{\hskip\hsize\hskip1em
\vbox{\hsize3.5cm\small\raggedright\pretolerance10000
 \noindent #1\hfill}\hss}\vbox to8pt{\vfil}\vss}}}

\begin{document}

\title{Uniqueness of curvature measures in pseudo-Riemannian geometry}

\author{Andreas Bernig}
\author{Dmitry Faifman}
\author{Gil Solanes}

\email{bernig@math.uni-frankfurt.de}
\email{faifmand@tauex.tau.ac.il} 
\email{solanes@mat.uab.cat}
\address{Institut f\"ur Mathematik, Goethe-Universit\"at Frankfurt,
Robert-Mayer-Str. 10, 60629 Frankfurt, Germany}
\address{School of Mathematical Sciences, Tel Aviv University, Tel Aviv 6997801, Israel}
\address{Departament de Matem\`atiques, Universitat Aut\`onoma de Barcelona, 08193 Bellaterra, Spain}


\begin{abstract}
The recently introduced Lipschitz-Killing curvature measures on pseudo-Riemannian manifolds satisfy a Weyl principle, i.e. are invariant under isometric embeddings. We show that they are uniquely characterized by this property. We apply this characterization to prove a K\"unneth-type formula for Lipschitz-Killing curvature measures, and to classify the invariant generalized valuations and curvature measures on all isotropic pseudo-Riemannian space forms.
\end{abstract}

\thanks{{\it MSC classification}:  53C65, 
53C50 
\\ A.B. was supported
 by DFG grant BE 2484/5-2\\
 G.S. was supported by FEDER/MICINN grant PGC2018-095998-B-I00 and the Serra H\'unter Programme}
\maketitle

\tableofcontents

\section{Introduction}

\subsection{Background}

A valuation on a finite-dimensional vector space $V$ is a functional $\mu:\mathcal{K}(V) \to A$, where $\mathcal{K}(V)$ denotes the set of compact convex subsets of $V$ and $A$ is an abelian semigroup, such that 
\begin{displaymath}
\mu(K \cup L)+\mu(K \cap L)=\mu(K)+\mu(L)
\end{displaymath}
whenever $K,L,K \cup L \in \mathcal{K}(V)$. 

An important example is given by the intrinsic volumes. If $V$ is a Euclidean vector space of dimension $n$, $K \in \mathcal K(V)$, then, as observed by Steiner \cite{steiner}, the volume of the $r$-tube $K_r:=K+rB$ around $K$ is a polynomial in $r$:
\begin{displaymath}
\vol K_r=\sum_{k=0}^n \mu_k(K)\omega_{n-k}r^{n-k}.
\end{displaymath}
Here $\omega_{n-k}$ is the volume of the $(n-k)$-dimensional unit ball. The coefficient $\mu_k(K)$ is called \emph{$k$-th intrinsic volume}. If $\iota:V \to W$ is an isometric embedding of Euclidean vector spaces, then $\mu_k^W(\iota(K))=\mu_k^V(K)$ for all $K \in \mathcal K(V)$. In particular, $\mu_k$ is invariant under translations and rotations. Conversely, if $\mu$ is a continuous (with respect to the Hausdorff metric on $\mathcal K(V)$), translation- and rotation-invariant real-valued valuation, then $\mu$ is a linear combination of intrinsic volumes by a famous theorem of Hadwiger.

Hadwiger's theorem has inspired a lot of research. To mention just a few of the numerous results, we refer the reader to \cite{ alesker03_un, bernig_sun09, bernig_g2, bernig_qig, bernig_solanes, bernig_voide} for versions for subgroups of the orthogonal group, to \cite{alesker_bernig_schuster, hug_schneider_schuster_a, ludwig_2005, schuster10, wannerer_area_measures, wannerer_unitary_module} for valuations taking values in some abelian semigroups, and to \cite{ludwig_reitzner99, ludwig_reitzner10} for semi-continuous valuations. 

A differential geometric version of Steiner's formula was found by Weyl  \cite{weyl_tubes}. Instead of taking a compact convex body, he considered a compact submanifold $M$ (possibly with boundary) of a Euclidean space and showed that the volume of an $r$-tube is a polynomial for small enough $r$. Moreover, the coefficients only depend on the intrinsic geometry of the submanifold, and not on the embedding. We refer to this as \emph{Weyl's principle}. 

For both formulas, the Steiner and the Weyl formula, local versions exist, where one looks only at those points in the $r$-tube such that the foot point on $K$ or $M$ belongs to a given Borel subset of $V$. The coefficients $\Lambda_k, k=0,\ldots,n$ are then valuations with values in the space of signed measures on $V$ and are called \emph{Lipschitz-Killing curvature measures}. For instance, if $(M,g)$ is a compact $m$-dimensional Riemannian submanifold without boundary, then $\Lambda_{m-2}^M(M,U)=\frac{1}{4\pi} \int_U  {sc} \cdot d\mathrm{vol}$, where $U \subset M$ is a Borel subset and $ {sc}$ is the scalar curvature of $(M,g)$. 

The structural similarity of the results by Steiner and Weyl is not a coincidence and can be explained with Alesker's much more recent theory of \emph{valuations on manifolds}. In this language, the intrinsic volumes are valuations which are defined on arbitrary Riemannian manifolds and which behave naturally with respect to isometric embeddings; and the Lipschitz-Killing curvature measures are curvature measures naturally associated to Riemannian manifolds. Conversely, Fu and Wannerer \cite{fu_wannerer} have recently shown in the spirit of Hadwiger's characterization that the intrinsic volumes/Lipschitz-Killing curvature measures are characterized by the Weyl principle, i.e. any valuation/curvature measure on Riemannian manifolds that satisfies the Weyl principle is a linear combination of intrinsic volumes/Lipschitz-Killing measures.  

It is a very natural question to look for analogous results in pseudo-Riemannian geometry. In this case, the tubes are in general not compact, but nevertheless one can try to associate valuations and curvature measures to pseudo-Riemannian manifolds. In the flat case, this was achieved in \cite{alesker_faifman} and \cite{bernig_faifman_opq}. It turns out that the continuity assumption is too restrictive and should be replaced by the notion of \emph{generalized} valuations or curvature measures. The very rough idea is that a generalized valuation can not be evaluated on every compact differentiable polyhedron, but only on smooth enough sets which are transversal (in some precise sense) to the light cone of the metric. These sets are called \emph{LC-transversal}. One obtains a sequence $\mu_k, {k\geq 0}$ of  complex-valued generalized valuations called \emph{intrinsic volumes}, and a sequence $\Lambda_k, {k\geq0}$ of  complex-valued generalized curvature measures called \emph{Lipschitz-Killing curvature measures}. 

The extension to the curved case was carried out in \cite{bernig_faifman_solanes}. To every pseudo-Riemannian manifold $M^{p,q}$ we associate a space $\mathcal{LK}(M)$ of \emph{intrinsic volumes}, which are generalized valuations on $M$; and a space $\widetilde{\mathcal{LK}}(M)$ of \emph{Lipschitz-Killing curvature measures}, which are generalized curvature measures on $M$. They can be evaluated on smooth enough LC-transversal sets. The most important property of these objects is the Weyl principle, which states that for every isometric immersion $M \looparrowright N$ of pseudo-Riemannian manifolds, the restriction of the intrinsic volume $\mu_k^N$ to $M$ equals $\mu_k^M$ and the restriction of the Lipschitz-Killing curvature measure $\Lambda_k^N$ to $M$ equals $\Lambda_k^M$. A characterization of the intrinsic volumes as the only generalized valuations on pseudo-Riemannian manifolds satisfying a Weyl principle appeared in \cite[Theorem D]{bernig_faifman_solanes}, based on the results from \cite{bernig_faifman_opq}.

\subsection{Results}

To complete the analogy with the Euclidean/Riemannian case, we still need a characterization theorem for generalized curvature measures satisfying a Weyl principle. 
As noted above, each Lipschitz-Killing curvature measure satisfies the Weyl principle, i.e. it associates to each pseudo-Riemannian manifold $(M,Q)$ a generalized curvature measure $\Lambda_k^M$ such that whenever $M \looparrowright N$ is an isometric immersion, then $\Lambda_k^N|_M=\Lambda_k^M$. The same holds true for linear combinations $\sum_{k=0}^\infty a_k \Lambda_k+b_k\bar \Lambda_k$. Our first main theorem states that, conversely, every assignment of a generalized curvature measure $\Lambda^M$ to each pseudo-Riemannian manifold $M$ that satisfies $\Lambda^N|_M=\Lambda^M$ is of this form. In other words, the space $\widetilde{\mathcal{LK}}(M)$ of Lipschitz-Killing curvature measures is characterized by the Weyl principle.

To state the theorem more precisely, we need some terminology. Let $\mathbf{\Psi Met}$ denote the category of pseudo-Riemannian manifolds with isometric immersions. Let $\mathbf{GCrv}$ be the category where the objects are pairs $(M, \Phi)$, with $M$ a smooth manifold, and $\Phi \in \mathcal C^{-\infty}(M, \mathbb C)$ (the space of generalized curvature measures). The morphisms $e:(M, \Phi_M)\to (N,\Phi_N)$ are immersions $e:M \looparrowright N$ such that $e^*\Phi_N$ is well-defined, and $\Phi_M=e^*\Phi_N$. The category $\mathbf{GVal}$ of manifolds with generalized valuations is defined similarly.

 A \emph{Weyl functor} on pseudo-Riemannian manifolds with values in generalized curvature measures is any covariant functor $\Lambda:\mathbf{\Psi Met}\to  \mathbf {GCrv}$ intertwining the forgetful functor to the category of smooth manifolds. More generally, we may similarly define Weyl functors between any two categories of manifolds equipped with a geometric structure, when natural restriction operations are  available for both structures. Important examples of Weyl functors are the intrinsic volumes of Riemannian manifolds, taking values in smooth valuations, and the Lipschitz-Killing curvature measures \cite{federer59, fu_wannerer}. For a different example, a family of Weyl functors on contact manifolds with values in generalized valuations was described in \cite{faifman_contact}.

 In this language, the intrinsic volumes and Lipschitz-Killing curvature measures on Riemannian manifolds were extended in \cite{bernig_faifman_solanes} to Weyl functors $\mu_k:\mathbf{\Psi Met}\to \mathbf{GVal}$ and $\Lambda_k:\mathbf{\Psi Met}\to \mathbf{GCrv}$, respectively. It was moreover shown that the Weyl functors $\mathbf{\Psi Met}\to \mathbf{GVal}$ are spanned over $\R$ by $\mu_0=\chi$ and $\{\mu_k, \overline{\mu_k}\}_{k\geq 1}$. Our first result is a similar classification for the curvature measures.

\begin{MainTheorem} \label{mainthm_uniqueness_functors}
Any Weyl functor $\Lambda:\mathbf{\Psi Met} \to  \mathbf {GCrv}$ is given by a unique infinite linear combination $\Lambda=\sum_{k=0}^\infty a_k \Lambda_k+b_k\bar \Lambda_k$.
\end{MainTheorem}

Theorem \ref{mainthm_uniqueness_functors} may be used to prove geometric formulas by the \emph{template method}, where the templates are special pseudo-Riemannian manifolds. As a first application of this method, we prove the following formula, which is well-known in the Riemannian case \cite[Equation 3.34]{cheeger_mueller_schrader}.

\begin{MainTheorem}[K\"unneth-type formula] \label{mainthm_kuenneth}
 Let $(M_1,Q_1),(M_2,Q_2)$ be pseudo-Riemannian manifolds. Let $A_i \subset M_i, i=1,2$ be LC-transversal differentiable polyhedra. 
 Then 
\begin{displaymath}
\Lambda_k^{M_1 \times M_2}(A_1\times A_2,\bullet)=\sum_{k_1+k_2=k} \Lambda_{k_1}^{M_1}(A_1,\bullet) \boxtimes \Lambda_{k_2}^{M_2}(A_2,\bullet).
\end{displaymath}
Here $\Lambda_{k_1}^{M_1}(A_1,\bullet) \boxtimes \Lambda_{k_2}^{M_2}(A_2,\bullet)$ denotes the exterior product of the generalized measures $\Lambda^{M_1}_{k_1}(A_1,\bullet)$ and $\Lambda^{M_2}_{k_2}(A_2,\bullet)$. 
\end{MainTheorem}

By the Weyl principle, we have for every pseudo-Riemannian manifold 
\begin{displaymath}
\mathcal{LK}(M) \subset \mathcal{V}^{-\infty}(M)^{\mathrm{Isom}(M)}, \quad \widetilde{\mathcal{LK}}(M) \subset \mathcal{C}^{-\infty}(M)^{\mathrm{Isom}(M)},
\end{displaymath}
 where $\mathrm{Isom}(M)$ is the isometry group.

A (pseudo-Riemannian) space form is a complete connected pseudo-Riemannian manifold of constant sectional curvature. We refer to Section \ref{subsec_space_forms} for the classification of space forms. A connected space form is \emph{isotropic} if the isometry group acts transitively on the level sets of the metric. Examples include all pseudospheres, pseudohyperbolic spaces, and flat pseudo-Euclidean spaces. Our third main theorem states that the displayed inclusions become equalities if $M$ is an isotropic space form. It thus gives a complete description of isometry invariant generalized valuations and curvature measures for these space forms, and can be considered as a Hadwiger-type theorem. 

\begin{MainTheorem} \label{mainthm:classification_spq}
Let $(M,Q)$ be an isotropic space form. Then 
\begin{align*}
\mathcal{V}^{-\infty}(M)^{\Isom(M)} & =\mathcal{LK}(M),\\
\mathcal{C}^{-\infty}(M)^{\Isom(M)} & =\widetilde{\mathcal{LK}}(M). 
\end{align*}
\end{MainTheorem}

The following special case of the theorem completes the classification of isometry invariant generalized valuations from \cite{bernig_faifman_opq} and will be the main ingredient in the proof of Theorem \ref{mainthm:classification_spq}.

\begin{Proposition}[Classification of isometry invariant curvature measures on pseudo-Euclidean space] \label{prop_flat_case}
Let $V$ be an $n$-dimensional real vector space endowed with a non-degenerate quadratic form of signature $(p,q)$ and isometry group $\OO(p,q)$. Let $\Curv_k^{-\infty}(V)^{\OO(p,q)}$ be the space of $k$-homogeneous, translation-invariant and $\OO(p,q)$-invariant generalized curvature measures on $V$. Then
\begin{equation} \label{eq_dimension_formula}
\dim\Curv^{-\infty}_k(V)^{\OO(p,q)} = \begin{cases} 2 & \text{ if } p,q \geq 1 \text{ and }  0 \leq k \leq n-1,\\ 1 & \text{ if } \min(p,q)=0 \text{ or } k=n. \end{cases}
\end{equation}
In each case, a basis is given by the real and/or imaginary parts of the Lipschitz-Killing curvature measure $\Lambda_k$.
\end{Proposition}

\subsection{Acknowledgements}
Part of this work was carried out during the second named author's stay at CRM - Universit\'e de Montr\'eal, which is gratefully acknowledged.
\section{Preliminaries}

\subsection{Pseudo-Riemannian space forms}
\label{subsec_space_forms}

We refer to \cite{oneill,wolf61} for the material in this subsection.
\begin{Definition}
\begin{enumerate}
\item The \emph{pseudo-Euclidean space of signature $(p,q)$} is $\R^{p,q}=\R^{p+q}$ with a quadratic form of signature $(p,q)$, e.g. $Q=\sum_{i=1}^p dx_i^2-\sum_{i=p+1}^{p+q} dx_i^2$. 
\item The \emph{pseudosphere of signature $(p,q)$ and radius $r>0$} is 
\begin{displaymath}
S_r^{p,q}=\{v \in \R^{p+1,q}:Q(v)=r^2\}.
\end{displaymath}
The pseudosphere $S_1^{n,1}\subset \R^{n+1,1}$ is called \emph{de Sitter space} and denoted by $dS^{n,1}$. 
\item The \emph{pseudohyperbolic space of signature $(p,q)$ and radius $r>0$} is 
\begin{displaymath}
H_r^{p,q}=\{v \in \R^{p,q+1}:Q(v)=-r^2\}.
\end{displaymath}
The pseudohyperbolic space $H_1^{n,1}$ is called the \emph{anti-de Sitter space}.
\end{enumerate}
\end{Definition}

The isometry groups of these spaces are given by 
\begin{align*}
\Isom(\R^{p,q}) & \cong \overline{\OO}(p,q)=\OO(p,q) \ltimes \R^{p,q},\\
\Isom(S_r^{p,q}) & \cong \OO(p+1,q),\\
\Isom(H_r^{p,q}) & \cong \OO(p,q+1).
\end{align*}
In each case, the action is transitive and the stabilizer at any point is conjugate to $\OO(p,q)$.

\begin{Definition}
A complete connected pseudo-Riemannian manifold of constant sectional curvature is called \emph{space form}. A connected space form whose isometry group acts transitively on the level sets of the metric is called \emph{isotropic}.
\end{Definition}

By a theorem of Wolf \cite{wolf61}, the stabilizer of a point in the isometry group of an isotropic space form acts on the tangent space by the full orthogonal group.

The next theorem gives a classification of simply connected space forms. They are all isotropic. 

\begin{Theorem} \label{thm_universal_cover_space_form}
Let $(M,Q)$ be a pseudo-Riemannian space form of signature $(p,q)$ and curvature $K$. Then the universal cover of $M$ is isometric to 
\begin{enumerate}
\item $\R^{p,q}$ if $K=0$;
\item $S^{p,q}_{\frac{1}{\sqrt K}}$ if $K>0$ and $p \geq 2$;
\item $\tilde S^{1,q}_{\frac{1}{\sqrt K}}$ (simply connected pseudo-Riemannian covering of $S^{1,q}_{\frac{1}{K}}$) if $K>0, p=1$;
\item $cS^{0,q}_{\frac{1}{\sqrt K}}$ (connected component of $S^{0,q}_{\frac{1}{\sqrt K}}$) if $K>0, p=0$;
\item $H^{p,q}_{\frac{1}{\sqrt{-K}}}$ if $K<0,q \geq 2$;
\item $\tilde H^{p,1}_{\frac{1}{\sqrt{-K}}}$ (simply connected pseudo-Riemannian covering of $H^{p,1}_{\frac{1}{\sqrt{-K}}}$) if $K<0,q=1$;
\item $cH^{p,0}_{\frac{1}{\sqrt{-K}}}$ (connected component of $H^{p,0}_{\frac{1}{\sqrt{-K}}}$) if $K<0, q=0$.
\end{enumerate}
\end{Theorem}

\subsection{Valuations and curvature measures on manifolds}

In this subsection, we recall the definitions of the basic objects of this paper: smooth and generalized valuations and smooth and generalized curvature measures on manifolds. We refer to \cite[Section 2]{bernig_faifman_solanes} for more details.

Let $M$ be a smooth manifold, assumed oriented for simplicity. By $\pi:\p_M \to M$ we denote the cosphere bundle. We let $\mathcal P(M)$ denote the set of compact differentiable polyhedra. A \emph{smooth valuation on a manifold} $M$ is a functional $\mu:\mathcal P(M) \to \R$ of the form 
\begin{equation} \label{eq_def_smooth_val}
 \mu(A)=\int_A \phi+\int_{\nc(A)} \omega, \quad \phi \in \Omega^n(M), \omega \in \Omega^{n-1}(\p_M).
\end{equation}

Here $\nc(A)$ is the normal cycle of $A$, which is an integral current in $\p_M$. The space of smooth valuations is a Fr\'echet space which is denoted by $\mathcal V^\infty(M)$. Examples are the Euler characteristic $\chi$, the volume, and more generally the intrinsic volumes of a Riemannian manifold $(M,g)$.

There is a natural notion of the emph{support} of a valuation, and the space of compactly supported smooth valuations is denoted by $\mathcal V^\infty_c(M)$ and equipped with a natural LF-topology. A \emph{generalized valuation} is an element of the dual space $\mathcal V^{-\infty}(M):=\mathcal{V}^\infty_c(M)^*$. It can be represented by generalized forms $\phi \in \Omega^n_{-\infty}(M), \omega \in \Omega^{n-1}_{-\infty}(\p_M)$. In this case, \eqref{eq_def_smooth_val} still makes sense for particularly nice $A$. Every $A \in \mathcal P(M)$ defines a generalized valuation $\chi_A$ by setting $\langle \chi_A, \mu\rangle=\mu(A)$, and every smooth valuation can also be considered as a generalized valuation by Alesker-Poincar\'e duality. The \emph{wave front} of a generalized valuation describes its singularities, we refer to \cite[Section 8]{alesker_bernig} for the definition. Given closed subsets $\Lambda \subset \p_+(T^*M), \Gamma \subset \p_+(T^*\p_M)$, the space of generalized valuations with wave front included in $(\Lambda,\Gamma)$ is denoted by $\mathcal V^{-\infty}_{\Lambda,\Gamma}(M)$.

A \emph{smooth curvature measure} is a functional of the form 
\begin{equation} \label{eq_def_smooth_curv}
 \Phi(A,U)=\int_{A \cap U} \phi+\int_{\nc(A) \cap \pi^{-1}(U)} \omega, \quad \phi \in \Omega^n(M), \omega \in \Omega^{n-1}(\p_M).
\end{equation}
Here $A \in \mathcal P(M)$ and $U \subset M$ is a Borel subset. The Fr\'echet space of smooth curvature measures is denoted by $\mathcal C^\infty(M)$. Sometimes we also write $[\phi,\omega]$ for the curvature measure defined by \eqref{eq_def_smooth_curv}.

The pairs of forms $(\phi,\omega)$ such that the valuation $\mu$ from \eqref{eq_def_smooth_val} is trivial were described in \cite{bernig_broecker07} in terms of the contact structure on $\p_M$. We need a (simpler) version of this description for curvature measures. 

Note that a local contact form $\alpha$ is unique up to multiplication by a non-zero function. In the following, we will do some  constructions using $\alpha$, and will leave it to the reader to check that each construction is independent of the choice of $\alpha$. 

A form $\omega \in \Omega^k(\p_M)$ is called \emph{primitive} if $k \leq n-1$ and 
 \begin{displaymath}
  d\alpha^{n-k} \wedge \omega \equiv 0 \mod \alpha.
 \end{displaymath}

The \emph{Lefschetz decomposition} of a form $\omega \in \Omega^k(\p_M)$ is given by  
 \begin{displaymath}
  \omega \equiv \sum_{i=0}^{\left\lfloor \frac{k}{2}\right\rfloor} \omega_i \mod \alpha,
 \end{displaymath}
where $\omega_i \in \Omega^{k}(\p_M)$ is of the form $\omega_i=d\alpha^i \wedge \tilde \omega_i$ with $\tilde \omega_i \in \Omega^{k-2i}(\p_M)$ primitive. See \cite{huybrechts05} for these notions.

\begin{Proposition} \label{prop_kernel_smooth}
Let $\phi \in \Omega^n(M), \omega \in \Omega^{n-1}(\p_M)$. 
The following conditions are equivalent.
\begin{enumerate}
 \item The curvature measure $\Phi$ defined by \eqref{eq_def_smooth_curv} vanishes. 
 \item $\phi=0$ and $\omega$ belongs to the ideal generated by $\alpha$ and $d\alpha$. 
 \item $\phi=0$ and 
 \begin{displaymath}
 \int_{\p_M} \alpha \wedge \omega \wedge \tau_0=0 
 \end{displaymath}
for all primitive forms $\tau_0 \in \Omega_c^{n-1}(\p_M)$.
\end{enumerate}
\end{Proposition}

\proof 
\begin{enumerate}
\item[$ $] 
\item[ $ii) \implies i)$] This follows from the fact that normal cycles are Legendrian.
 \item[$i) \implies ii)$] Take a smooth compact submanifold $A \subset M$ of dimension $n$ with smooth boundary $\partial A$. Then for every $f \in C_c^\infty(M)$, we have 
\begin{displaymath}
\int_A f \phi+\int_{\nc(A)} \pi^*f \cdot \omega=0.
\end{displaymath}
Taking $f$ be supported in $\mathrm{int}A$, we find that $\phi=0$. Letting the support of $f$ shrink to a point on the boundary, it follows that $\omega$ vanishes on all tangent spaces to $\nc(A)$. Since these tangent spaces are dense in the set of all Legendrian planes, it follows that $\omega$ vanishes on all Legendrian planes. By \cite[Lemma 1.4]{bernig_broecker07} this implies that $\omega \in \langle \alpha,d\alpha\rangle$. 
 \item[$ii) \implies iii)$] Obvious.
 \item[$iii) \implies ii)$] Let $\tau \in \Omega_c^{n-1}(\p_M)$ be arbitrary and let 
\begin{displaymath}
  \omega \equiv \sum_{i=0}^{\left\lfloor \frac{n-1}{2}\right\rfloor} \omega_i \mod \alpha, \quad \tau \equiv \sum_{i=0}^{\left\lfloor \frac{n-1}{2}\right\rfloor} \tau_i \mod \alpha, 
 \end{displaymath}
be the Lefschetz decompositions. Then $\omega_i \wedge \tau_0 \equiv 0$ and $\tau_i \wedge \omega_0 \equiv 0$ for all $i>0$. 
The assumption is thus equivalent to  
\begin{displaymath}
0= \int_{\p_M} \alpha \wedge \omega \wedge \tau_0 = \int_{\p_M} \alpha \wedge \omega_0 \wedge \tau_0 = \int_{\p_M} \alpha \wedge \omega_0 \wedge \tau, 
\end{displaymath}
which implies by Poincar\'e duality that $\alpha \wedge \omega_0=0$. Hence $\omega_0$ is a multiple of $\alpha$. Since each $\omega_i, i>0$ is a multiple of $d\alpha$, the statement follows.
\end{enumerate}
\endproof

A \emph{generalized curvature measure} is given by a pair $\phi \in \Omega^n_{-\infty}(M), \omega \in \Omega^{n-1}_{-\infty}(\p_M)$. It can be evaluated at pairs $\mu,f$, where $\mu \in \mathcal V^\infty_c(M)$ and $f \in C_c^\infty(M)$. We write $\Phi(\mu,f)$ for this evaluation. The space of generalized curvature measures is denoted by $\mathcal C^{-\infty}(M)$.  As for generalized valuations, the singularities of a generalized curvature measure can be described by its wave front set \cite[Section 2.3]{bernig_faifman_solanes}. The set of generalized curvature measures with wave front set contained in $(\Lambda,\Gamma)$, where $\Lambda \subset \p_+(T^*M), \Gamma \subset \p_+(T^*\p_M)$ are closed subsets, is denoted by $\mathcal C^{-\infty}_{\Lambda,\Gamma}(M)$.

If $A \in \mathcal P(M)$ satisfies certain transversality conditions (which are given in terms of wave front sets), then 
\begin{displaymath}
 \Phi(A,f)=\int_{A} f \phi+\int_{\nc(A)} \pi^*f \cdot \omega, \quad f \in C_c^\infty(M)
\end{displaymath}
is well-defined.

\begin{Proposition} \label{prop_kernel_generalized}
Let $\phi \in \Omega^n_{-\infty}(M), \omega \in \Omega^{n-1}_{-\infty}(\p_M)$. The following conditions are equivalent.
\begin{enumerate}
 \item The generalized curvature measure $\Phi$ defined by $(\phi,\omega)$ vanishes. 
 \item $\phi=0$ and $\omega$ belongs to the ideal in $\Omega^*_{-\infty}(\p_M)$ generated by $\alpha$ and $d\alpha$. 
 \item $\phi=0$ and 
 \begin{displaymath}
 \int_{\p_M} \alpha \wedge \omega \wedge \tau_0=0 
 \end{displaymath}
for all primitive forms $\tau_0 \in \Omega_c^{n-1}(\p_M)$.
\end{enumerate}
\end{Proposition}

\proof
The implications $ii) \implies i), ii) \iff iii)$ are as in the proof of Proposition \ref{prop_kernel_smooth}.

For the implication $i) \implies ii)$, we work locally and with coordinates and assume that $M=\R^n$. Convolve with an approximate identity $\rho_j \in C_c^\infty(\mathrm{GL}(n))$. Then $\rho_j * \Phi$ is the smooth curvature measure represented by the smooth forms $(\rho_j \ast \phi, \rho_j * \omega)$, but obviously it is the trivial curvature measure. By Proposition \ref{prop_kernel_smooth}, $\rho_j * \phi=0$, while $\rho_j * \omega$ belongs to the ideal $\langle \alpha,d\alpha\rangle \subset \Omega^{n-1}(\p_M)$. For $j \to \infty$, $\rho_j * \phi \to \phi$ and $\rho_j * \omega \to \omega$ in the weak topology, hence $\phi=0$. Since $\langle \alpha,d\alpha\rangle \cap \Omega^{n-1}_{-\infty}(\p_M)$ is closed in the weak topology (which follows from the implication $ii) \iff iii)$), it follows that $\omega$ belongs to this space. 
\endproof

Sometimes we also need $\C$-valued valuations and curvature measures, which are defined in an analogous way. In all the following, the range which is either $\R$ or $\C$ is often omitted from notation, and should be determined from context.

\subsection{Translation-invariant valuations and curvature measures}

If $V$ is a vector space of dimension $n$, the spaces of smooth or generalized translation-invariant valuations and curvature measures are denoted by $\Val^{\pm \infty},\Curv^{\pm \infty}$. They admit gradings by homogeneity
\begin{align*}
\Val^{\pm \infty} & =\bigoplus_{k=0}^n \Val^{\pm \infty}_k,\\
\Curv^{\pm \infty} & =\bigoplus_{k=0}^n \Curv^{\pm \infty}_k.
\end{align*}

A $k$-homogeneous element in one of these spaces can be represented by a pair $(0,\omega)$ with $\omega$ translation-invariant and of bidegree $(k,n-k-1)$ if $k<n$; and by a pair $(\phi,0)$ with $\phi$ translation-invariant if $k=n$.

\begin{Proposition} \label{prop_kernel_translation_inv}
The (smooth or generalized) curvature measure induced by a translation-invariant (smooth or generalized) form $\omega$ of bidegree $(k,n-k-1)$ with $k<n$ vanishes if and only if $\omega$ belongs to the ideal in $\Omega^{*}(\p_V)^{\tr}$ or $\Omega^{*}_{-\infty}(\p_V)^{\tr}$ generated by $\alpha$ and $d\alpha$. 
\end{Proposition}

The proof is similar to the proofs of Propositions \ref{prop_kernel_smooth} and \ref{prop_kernel_generalized} and we omit the details. In this case, instead of taking the usual Poincar\'e pairing on manifolds, we have to use the Poincar\'e pairing on translation-invariant forms as follows. Take the wedge product of two translation-invariant forms of complementary degrees, and push-forward to $V$. Then we obtain a translation-invariant $n$-form on $V$, hence a multiple of the volume form. The corresponding factor is then the pairing of the two forms. This pairing is non-degenerate. 

\subsection{LC-transversality}

Let $(M,Q)$ be a pseudo-Riemannian manifold.  We denote by $\LC^*_M \subset \p_M=\p_+(T^*M)$ the set of null-directions in the cosphere bundle and by $\LC_M \subset \p_+(TM)$ the set of null-directions in the sphere bundle. For a submanifold $X \subset M$, we let $N^*X \subset T^*M$ be the conormal bundle, which we consider often as a subset of $\p_+(T^*M)$. Using the metric $Q$, we may identify this set with a subset in $TM$ (or $\p_+(TM)$), denoted by $N^QX$.

\begin{Definition}[{\cite[Section 4.2]{bernig_faifman_solanes}}]
A differentiable polyhedron $A\subset M$ is called \emph{LC-transversal} if each smooth stratum of $\nc(A)$ intersects $\LC^*_M$ transversally. In particular, a submanifold $X \subset M$ is LC-transversal if $N^QX \pitchfork \LC_M$ in $\mathbb P_+(TM)$.
\end{Definition}

We will need the following generalization of LC-transversality.

\begin{Definition}
A generalized valuation $\psi \in \mathcal{V}^{-\infty}(M)$ is called \emph{LC-transversal} if $\psi \in \mathcal V^{-\infty}_{\Lambda,\Gamma}(M)$ with $\Gamma \cap N^*(\LC_M^*)=\emptyset$. 
\end{Definition}

This relates to the LC-transversality of subsets as follows.
\begin{Lemma}\label{lem:lc_transversal_implies_lc_transversal}
	If a differentiable polyhedron $A \subset M$ is LC-transversal, then the generalized valuation $\chi_A$ is LC-transversal.
\end{Lemma}
\proof
Recall that $\chi_A=[([[A]], [[\nc(A)]]  )]$. Now $\WF([[\nc(A)]])$ is contained in the union of the conormal bundles to the smooth strata of $\nc(A)$. By assumption, the latter conormal bundles have empty intersection with $N^*\LC_M^*$.
\endproof

\section{Uniqueness of the Lipschitz-Killing functors}

In this section we will prove Proposition \ref{prop_flat_case}, which is the technical heart of the proof of Theorems \ref{mainthm_uniqueness_functors} and \ref{mainthm:classification_spq}. We will need two technical lemmas. The first one is from \cite{bernig_faifman_opq}, see also \cite[Section 4.4]{alesker_faifman} for some of the notation. To state it, we need some preparation. 

Let $X$ be a smooth manifold, and $E$ a smooth vector bundle over $X$.
For any $\nu \geq 0$ and a locally closed submanifold $Y\subset X$, define the vector bundle $F^\nu_Y$ over $Y$ with fiber

\begin{displaymath}
F^\nu_Y|_y=\Sym^\nu(N_yY) \otimes \Dens^*(N_yY) \otimes E|_y.                                                                                                                                                                                                                                                                                                                                                                                                                                                                                                                                                                                    \end{displaymath}

For a closed submanifold $Y\subset X$, let $\Gamma^{-\infty,\nu}_{Y}(X,E) \subset \Gamma^{-\infty}_{Y}(X,E)$ be the space of all generalized sections supported on $Y$ with differential order not greater that $\nu\geq 0$ in directions normal to $Y$. One then has a natural isomorphism 
\begin{displaymath}
\Gamma_Y^{-\infty,\nu}(X,E)/\Gamma_Y^{-\infty,\nu-1}(X,E)\cong\Gamma^{-\infty}(Y,F_Y^\nu).
\end{displaymath}

Now let a Lie group $G$ act on $X$ in such a way that there are finitely many orbits, all of which are locally closed submanifolds. We will assume that $E$ is a $G$-vector bundle. If $Y \subset X$ is a $G$-invariant locally closed submanifold, then $F^\nu_Y$ is naturally a $G$-bundle. If $Y$ is in fact a closed submanifold then $\Gamma^{-\infty,\nu}_Y(X,E)^G$, $\nu\geq 0$, form a filtration on $\Gamma^{-\infty}_Y(X,E)^G$.

\begin{Lemma}[{\cite[Lemma A.1]{bernig_faifman_opq}}]
\label{lem:LocallyClosedOrbits}
Let $Z\subset X$ be a closed $G$-invariant subset.
Decompose $Z=\bigcup_{j=1}^J Y_j$ where each $Y_j$ is a $G$-orbit and fix an element $y_j \in Y_j$. Then
\begin{displaymath}
 \dim \Gamma^{-\infty}_Z(X,E)^G \leq \sum_{\nu=0}^\infty \sum _{j=1}^J \dim \Gamma^\infty(Y_j,F^\nu_{Y_j})^{G}=\sum_{\nu=0}^\infty \sum _{j=1}^J \dim \left(F^{\nu}_{Y_j}|_{y_j}\right)^{\Stab(y_j)}.
\end{displaymath}
\end{Lemma}

For a $G$-module $X$ and a character $\chi$ on $G$, we write $X^{G,\chi}=\{\omega \in X: g\omega=\chi(g)\omega , g\in G \}$ and call its elements $(G,\chi)$-invariant. By tensorizing all representations with the one-dimensional representation $\C$ on which $G$ acts by $\chi(g)$, a similar upper bound holds for $(G,\chi)$ resp. $(\Stab(y_j),\chi)$-invariant subspaces. We will use the character $\det:\mathrm{O}(p,q) \to \R$.

The second technical lemma concerns  generalized functions on the real line. Assume $\sigma\in C^\infty(-\epsilon,\epsilon)$, $\sigma(0)=0$ and $\sigma'(0)\neq 0$. Let $g_\alpha$, $\alpha\geq 0$ be a smooth family of smooth injective maps from $(-\epsilon,\epsilon)$ to itself, given by  $g_\alpha(x)=x b_\alpha(x)$ where $b_0(x)=1$, $b_\alpha(x)\leq 1$, $b_\alpha(x)$ is smooth in both variables, and $\left.\frac{d}{d \alpha}\right|_0 b_\alpha(0)\neq 0$. Define $\psi_\alpha(x)=\frac{\sigma(g_\alpha x)}{\sigma(x)}$ for $x\neq 0$, and $\psi_\alpha(0) =\lim_{x \to 0} \psi_\alpha(x)=g_\alpha'(0)=b_\alpha(0)$.

\begin{Proposition}
	\label{prop:uniqueness_of_homogeneous_extension}
Let $W$ be the space of generalized functions $f\in C^{-\infty}(-\epsilon,\epsilon)$ satisfying the equation $g_\alpha^*f=\psi_{\alpha}^{-m}\cdot f$ for all $|\alpha| \leq\alpha_0$, for some $\alpha_0>0$ and $m \in \mathbb N$. Then the subspace of $W$ of functions supported at $x=0$ is at most one-dimensional. Moreover, there is $\delta>0$ such that the space of restrictions of $W$ to $(-\delta,\delta)\setminus\{0\}$ is also at most one-dimensional.
\end{Proposition}

\proof
Write $\sigma(x)=x s(x)$ and $w(x)=\left.\frac{\partial}{\partial \alpha}\right|_0b_\alpha(x)$. Assume $w(x), s(x)\neq 0$ for $|x|\leq \epsilon_1\leq\epsilon$. We have 
\begin{displaymath}
\psi_\alpha(x)=\frac{g_\alpha x}{x}\frac{s(g_\alpha x)}{s(x)}=b_\alpha(x)\frac{s(g_\alpha x)}{s(x)}\Rightarrow
\left.\frac{\partial}{\partial \alpha}\right|_0\psi_\alpha(x)=w(x) +\frac{1}{s(x)}xs'(x)w(x). 
\end{displaymath}

Differentiating the functional equation at $\alpha=0$, we find
\begin{displaymath}
xw(x)f'=-mf w(x)\left(1+x\frac{s'(x)}{s(x)}\right) \iff xs(x)f'=-m(s(x)+xs'(x))f.	
\end{displaymath}

Write $F(\sigma)=f(x)$ in an interval $|x|<\epsilon_2\leq \epsilon_1$ where $\sigma$ is invertible. Then $f'(x)=F'(\sigma)(s(x)+xs'(x))$, and the equation becomes
\begin{displaymath}
\sigma  F'(\sigma)=-mF(\sigma),
\end{displaymath}
which is the equation of a $(-m)$-homogeneous function in a neighborhood of the origin, and the statement follows from \cite[Chapter 1, \S 3]{gelfand_shilov}. 
\endproof

\proof[Proof of Proposition \ref{prop_flat_case}]

A translation-invariant generalized curvature measure of degree $n$ is induced by a translation-invariant generalized $n$-form on $V$. Since translations act transitively on $V$, such forms are actually smooth and hence multiples of the volume form. Hence $\Curv_n^{-\infty,\OO(p,q)}$ is spanned by the volume. 

Let us assume in the following that $0 \leq k \leq n-1$. Consider first the case $\min(p,q)=0$, i.e. the case of a Euclidean (or anti-Euclidean) vector space. The action of the Euclidean motion group $\overline{\OO(p,q)}$ on the cosphere bundle is transitive, which implies that all invariant generalized forms are smooth. However, from the classification of the invariant smooth forms in \cite{fu90, fu98} it follows immediately that $\Curv^{-\infty,\OO(p,q)}_k \cong \Curv^{\infty,\OO(p,q)}_k \cong \Val_k^{\infty,\OO(p,q)}$, and the latter space is one-dimensional by Hadwiger's theorem. 

In the remaining case $\min(p,q)>0$ we proceed as in \cite[Section 5.3]{bernig_faifman_opq}. 

Let $\p_+(V^*):=V^*\setminus \{0\}/\R_+$ and let $\p_V:=V \times \p_+(V^*)$ be the cosphere bundle over $V$. For $\xi \in \p_+(V^*)$ we denote by $\xi^Q$ the $Q$-orthogonal complement, which is a hyperplane in $V$. Denote by $D^{k,l}$ the vector bundle over $\p_+(V^*)$ whose fiber over $\xi$ is given by 
\begin{displaymath}
D^{k,l}|_\xi=\largewedge^k(V^*) \otimes \largewedge^l(\xi^Q) \otimes \xi^l.
\end{displaymath} 
Then the space of translation-invariant forms is given by 
\begin{displaymath}
\Omega_{-\infty}^{k,l}(\p_V)^{\tr} \cong \Gamma^{-\infty}(\p_+(V^*),D^{k,l}).
\end{displaymath}

The space of vertical translation-invariant generalized forms is 
\begin{displaymath}
\Omega_{v,-\infty}^{k,l}(\p_V)^{tr} \cong \Gamma^{-\infty}(\p_+(V^*),D_v^{k,l}),
\end{displaymath}
where $D_v^{k,l}$ is the vector bundle with fiber 
\begin{displaymath}
D_v^{k,l}|_\xi=\xi \otimes \largewedge^{k-1}(V^*/\xi) \otimes \largewedge^l(\xi^Q) \otimes \xi^l.
\end{displaymath} 

The space of horizontal translation-invariant generalized forms (i.e. the quotient of all generalized translation-invariant forms by the vertical translation-invariant generalized forms) is 
\begin{displaymath}
 \Omega_{h,-\infty}^{k,l}(\p_V)^{tr} \cong \Gamma^{-\infty}(\p_+(V^*),D_h^{k,l}),
\end{displaymath}
where $D_h^{k,l}$ is the vector bundle with fiber 
\begin{displaymath}
D_h^{k,l}|_\xi=D^{k,l}/D_v^{k,l}=\largewedge^k(\xi^{Q*}) \otimes \largewedge^l(\xi^Q) \otimes \xi^l.
\end{displaymath} 

Let $D_p^{k,l}|_\xi \subset D_h^{k,l}|_\xi$ by the subspace of primitive elements. By Proposition \ref{prop_kernel_translation_inv} 
\begin{displaymath}
\Curv_k^{-\infty} \cong \Omega_{p,-\infty}^{k,n-k-1}(\p_V)^{tr} \cong \Gamma^{-\infty}(\p_+(V^*),D_p^{k,n-k-1}).
\end{displaymath}

Let $G:=\OO(p,q)$. Since the orientation of the conormal cycle depends on the choice of an orientation on $V$, we have 
\begin{displaymath}
\Curv_k^{-\infty,G} \cong \Gamma^{-\infty}(\p_+(V^*),D_p^{k,n-k-1})^{G,\det}.
\end{displaymath}

Recall from \cite[Proposition 4.9]{bernig_faifman_opq} that the action of $G$ on $\p_+(V^*)$ has two open orbits $M^+:=\{\xi:Q(\xi,\xi)>0\}, M^-:=\{\xi: Q(\xi,\xi)<0\}$ and one closed orbit $M^0:=\{\xi:Q(\xi,\xi)=0\}$, which is called \emph{light cone}. 

We claim that the space $\Gamma^{-\infty}_{M^0}(\p_+(V^*),D_p^{k,n-k-1})^{G,\det}$ of $(G,\det)$-invariant and primitive generalized forms supported on the light cone is trivial if $(n-k)$ is odd, and at most one-dimensional if $(n-k)$ is even. 

To prove the claim, we fix $\xi \in M^0$ and denote by $H \subset G$ the stabilizer of $\xi$. 
By Lemma \ref{lem:LocallyClosedOrbits} and \cite[Lemma 5.7]{bernig_faifman_opq},
\begin{align}
  \dim \Gamma^{-\infty}_{M^0} & (\p_+(V^*),D_p^{k,n-k-1})^{G,\det} \nonumber \\
  & \leq \sum_{\nu=0}^\infty \dim \left(\Sym^\nu N_\xi M^0 \otimes \Dens^*(N_\xi M^0) \otimes D_p^{k,n-k-1}|_\xi\right)^{H,\det} \nonumber\\
  & = \sum_{\nu=0}^\infty \dim \left(D_p^{k,n-k-1}|_\xi \otimes \xi^{-2\nu-2}\right)^{H,\det} \nonumber\\
  & \leq  \sum_{\nu=0}^\infty \dim \left(D_h^{k,n-k-1}|_\xi \otimes \xi^{-2\nu-2}\right)^{H,\det}. \label{eq_inequality_from_lemma}
\end{align}

Set 
\begin{displaymath}
U:=D_h^{k,n-k-1}|_\xi \otimes \xi^{-2\nu-2} \cong \largewedge^k (\xi^Q)^* \otimes \largewedge^{n-k-1} \xi^Q \otimes \xi^\beta
\end{displaymath}
with $\beta=n-k-2\nu-3$.  

Since we have an exact sequence 
\begin{displaymath}
0 \to \xi \to \xi^Q \to \xi^Q/\xi \to 0,
\end{displaymath}
we get an exact sequence 
\begin{displaymath}
0 \to (\xi^Q/\xi)^* \to (\xi^Q)^* \to \xi^* \to 0.
\end{displaymath}

Note that $(\xi^Q/\xi)^* \cong \xi^Q/\xi$, since the restriction of $Q$ to $\xi^Q/\xi$ is non-degenerate. The action of $H$ on this space is that of $\OO(p-1,q-1)$.

Hence 

\begin{displaymath}
U_u:=\largewedge^k (\xi^Q/\xi) \otimes \largewedge^{n-k-1} \xi^Q \otimes \xi^\beta \subset U.
\end{displaymath}

By \cite[Lemma 2.1]{bernig_faifman_opq}  we have
\begin{displaymath}
W_u:=U/U_u \cong \largewedge^{k-1} (\xi^Q/\xi) \otimes \largewedge^{n-k-1} \xi^Q \otimes \xi^{\beta-1}.
\end{displaymath} 

As in \cite{bernig_faifman_opq} we define the subspaces 
\begin{displaymath}
A_{k,l,\beta} := \largewedge^k(\xi^Q/\xi) \otimes \largewedge^l(\xi^Q/\xi) \otimes \xi^\beta.
\end{displaymath}

By \cite[Lemma 5.1]{bernig_faifman_opq} we have 
\begin{displaymath}
\dim A_{k,l,\beta}^{H,\det}=\begin{cases}
                                       0 & \text{ if } k+l \neq n-2 \text{ or } \beta \neq 0,\\
                                       1 & \text{ if } k+l=n-2 \text{ and } \beta=0.
                                      \end{cases}
\end{displaymath}

Note that the kernel of $\largewedge^k \xi^Q \to \largewedge^k(\xi^Q/\xi)$ is canonically isomorphic to $\largewedge^{k-1}(\xi^Q/\xi)\otimes\xi$. Define the subspaces
\begin{align*}
U_{uu} & :=\largewedge^k (\xi^Q/\xi) \otimes \largewedge^{n-k-2} (\xi^{Q}/\xi) \otimes \xi^{\beta+1} \cong A_{k,n-k-2,\beta+1} \subset U_u, \\
U_{wu} & :=\largewedge^{k-1} \xi^Q/\xi \otimes \largewedge^{n-k-2} (\xi^Q/\xi) \otimes \xi^{\beta} \cong A_{k-1,n-k-2,\beta} \subset W_u,
\end{align*}
and the quotients 
\begin{align*}
W_{uu} & :=U_u/U_{uu} \cong \largewedge^k (\xi^Q/\xi) \otimes \largewedge^{n-k-1} (\xi^{Q}/\xi) \otimes \xi^\beta \cong A_{k,n-k-1,\beta}, \\
W_{wu} & := W_u/U_{wu} \cong \largewedge^{k-1} (\xi^Q/\xi) \otimes \largewedge^{n-k-1} (\xi^Q/\xi) \otimes \xi^{\beta-1} \cong A_{k-1,n-k-1,\beta-1}.
\end{align*}

If $A$ is an $H$-representation and $B$ a subrepresentation, then $\dim A^{H,\det} \leq \dim B^{H,\det} + \dim(A/B)^{H,\det}$. We thus obtain that 
\begin{align*}
 \dim W_u^{H,\det} & \leq \dim U_{wu}^{H,\det}+\dim W_{wu}^{H,\det} \\
 & = \dim A_{k-1,n-k-2,\beta}^{H,\det}+\dim A_{k-1,n-k-1,\beta-1}^{H,\det} \\
 & = \begin{cases} 0 & \beta \neq 1\\ 1 & \beta=1.\end{cases}
\end{align*}

Hence $W_u$ can only contain an $(H,\det)$-invariant if $\beta=1$. Let us show by a direct argument as in \cite[Prop. 5.10]{bernig_faifman_opq} that also in the case $\beta=1$, there is no such invariant. 

Let $\rho \in W_u^{H,\det}$. Let $Y \subset \xi^Q$ be a complement of $\xi$ and $H_Y \subset H:=\Stab(\xi)$ be the stabilizer of $Y$. The decomposition
\begin{displaymath}
 W_{u,\beta=1}=(\largewedge^{k-1} \xi^Q/\xi \otimes \largewedge^{n-k-1} Y) \oplus (\largewedge^{k-1} \xi^Q/\xi \otimes \largewedge^{n-k-2} Y \otimes \xi)
\end{displaymath}
is compatible with the action of $H_Y$. The projection of $\rho$ to the second summand is $H_Y$-invariant. Since there are elements in $H_Y$ acting by the identity on $Y$ (and hence on $\xi^Q/\xi$) and by rescaling $\xi$, the second summand can not contain non-zero $(H_Y,\det)$-invariant elements. Hence $\rho$ must belong to the first summand. Since $Y$ was an arbitrary complement of $\xi$, our invariant must belong to the intersection
\begin{displaymath}
 \bigcap_Y \largewedge^{k-1} \xi^Q/\xi \otimes \largewedge^{n-k-1} Y.
\end{displaymath}

Fix some $Y$ that gives a minimal-length representation 
\begin{displaymath}
 \rho=\sum_{i=1}^m \eta_i\otimes y_i, \quad \eta_i \in \largewedge^{k-1} \xi^Q/\xi, y_i \in \largewedge^{n-k-1} Y.
\end{displaymath}
By the minimality of the representation, the $\eta_i$ are linearly independent.

The subgroup of $H$ acting by $\mathrm{Id}$ on $\xi^Q/\xi$ is transitive on all hyperplanes $Y$ complementing $\xi$. Acting by such an element $g$ on $\rho$, we get the equality 
\begin{displaymath}
\sum \eta_i \otimes (y_i-y'_i)=0, \text{ where } y'_i=g(y_i) \in \largewedge^{n-k-1}Y'. 
\end{displaymath}

It follows that $y_i=y'_i$ for all $i$, hence $y_i \in \bigcap_Y \largewedge^{n-k-1}Y$. It is elementary to prove that this intersection is trivial, and hence $y_i=0$ for all $i$ and therefore $\rho=0$.

On the other hand, we have 
\begin{align*}
\dim U_u^{H,\det} & \leq \dim U_{uu}^{H,\det}+\dim W_{uu}^{H,\det} \\
 & =\dim A_{k,n-k-2,\beta+1}^{H,\det}+\dim A_{k,n-k-1,\beta}^{H,\det} \\
 & = \begin{cases} 0 & \beta \neq -1\\ 1 & \beta=-1,\end{cases}
\end{align*} 
and therefore 
\begin{displaymath}
\dim U^{H,\det} \leq \dim U_u^{H,\det}+\dim W_u^{H,\det} \begin{cases} =0 & \beta \neq -1\\ \leq 1 & \beta=-1.\end{cases}
\end{displaymath}

We distinguish two cases. If $(n-k)$ is odd, then $\beta=n-k-2\nu-3 \neq -1$ for all $\nu$. By \eqref{eq_inequality_from_lemma}, the space of $(G,\det)$-invariant generalized primitive forms of degree $(k,n-k-1)$ supported on the light cone is trivial. Let $\phi \in \Omega_{p,-\infty}^{k,n-k-1}(\p_V)^{tr}$ be $(G,\det)$-invariant. The restriction of $\phi$ to each open orbit of $\p_V$ must be a multiple of the form $\phi_{k,r}^\pm$ constructed in \cite[Section 5.1]{bernig_faifman_solanes}. Hence there are constants $c_\pm$ such that $\phi=c_+ \phi_{k,r}^+$ on $S^+V$ and $\phi=c_- \phi_{k,r}^-$ on $S^-V$. It follows that the global form $\phi-c_+\phi_{k,r}^0-c_-\phi_{k,r}^1$ is supported on the light cone. Since this form is $(G,\det)$-invariant, it vanishes. 

Let $(n-k)$ be even. Let $\phi \in \Omega_{p,-\infty}^{k,n-k-1}(\p_V)^{tr}$ be $(G,\det)$-invariant. In the following, we identify $V=V^*$ using the quadratic form $Q$. By the above, we obtain that $\phi=c_+ \phi_{k,r}^+$ on $S^+V$ and $\phi=c_- \phi_{k,r}^-$ on $S^-V$ for some constants $c_\pm$. We will show that $c_+=c_-$. 

Fix a compatible Euclidean structure $P$ on $V$ (see \cite[Definition 2.7]{bernig_faifman_opq}). From \cite[Equation (60)]{bernig_faifman_solanes},
we have 
\begin{displaymath}
 \phi_{k,r}^+=\sigma_+^{-\frac{n-k}{2}}\rho_{k,r},\quad \phi_{k,r}^-=-\sigma_-^{-\frac{n-k}{2}}\rho_{k,r},
\end{displaymath}
where $\rho_{k,r}$ is a globally defined smooth form, and $\sigma_+,\sigma_-$ are the positive and negative parts of $\sigma=\sigma_+-\sigma_-=\frac{Q}{P}$.

Write $V=\R^{p,q}=\R^{1,1} \oplus \R^{p-1,q-1}$ and $L:=\p_+(\R^{1,1}) \subset \p_+(V)$, $\xi_0, \xi'_0\in L$ the degenerate lines. Take $g_\alpha \in \SO^+(Q)$ fixing $\R^{p-1,q-1}$ and acting by an $\alpha$-boost on $\R^{1,1}$, with $g_\alpha|_{\xi_0}=e^\alpha$, $g_\alpha|_{\xi'_0}=e^{-\alpha}$. Thus $\{g_\alpha:\alpha\in\R\} \simeq \SO^+(1,1)$. Let $\psi_\alpha(\xi):=\left.\frac{g_\alpha^*\sigma}{\sigma}\right|_\xi$.

We will use the standard Euclidean structure on $\R^{1,1}$ and introduce the polar angle $\theta$ for which $\theta(\xi_0)=\frac{\pi}{4}$. Consider the interval $L'$ with $0\leq \theta\leq \frac{\pi}{2}$. 
Then $\sigma=\cos2\theta$ is a coordinate in the interior of $L'$. With respect to this coordinate, we have $g_\alpha(\sigma)=\sigma b_\alpha(\sigma)$ with
\begin{displaymath}
b_\alpha(\sigma)=\frac1{\cosh{2\alpha}+\sinh(2\alpha)\sqrt{1-\sigma^2}}. 
\end{displaymath}

Since $\phi_{k,r}=\sigma^{-\frac{n-k}{2}}\rho_{k,r}$ is $\SO(p,q)$-invariant, we find that 
\begin{equation} \label{eq_invariance_rho}
g_\alpha^*\rho_{k,r}=\psi_\alpha^{\frac{n-k}{2}}\rho_{k,r}. 
\end{equation}

By $\OO(p,q)$-invariance it follows that the wave front set of $\phi$ is disjoint from $N^*L$. Letting $j_L:L\hookrightarrow \mathbb P_+(V)$ be the inclusion and denoting $D=j_L^* D^{k,n-1-k}$, we may therefore define $\phi_L=j_L^*\phi\in\Gamma^{-\infty}(L, D)$.

By definition, $D^{k,n-k-1}=\largewedge^k V\otimes \largewedge^{n-k-1}\xi^Q \otimes \xi^{n-k-1}$. Under the action of $\SO^+(1,1)$, we can decompose into equivariant summands 
\begin{align*}
 V & =\R^{p-1,q-1}\oplus \xi_0\oplus\xi_0' \\
 \xi^Q & =\R^{p-1,q-1}\oplus\tilde \xi,
\end{align*}
where $\tilde\xi\subset\R^{1,1}$. The first summand can be further written as a sum of lines on which $\SO^+(1,1)$ acts trivially. 

It follows that $D$ can be $\SO^+(1,1)$-equivariantly decomposed into the sum of line bundles, $D=\oplus_{i=1}^N D_j$. Let $\pi_j:D \to D_j$ be the projection. Since $\rho_{k,r}$ is a smooth and non-vanishing form, there exists $j_0$ such that $\pi_{j_0}(\rho_{k,r})$ is a smooth non-vanishing section of $D_{j_0}$ over a neighborhood of $\xi_0$, denoted $\tilde L_0$. It follows that 
\begin{displaymath}
 \pi_{j_0}(\phi_L)=f \pi_{j_0}(\rho_{k,r})
\end{displaymath}
for some $f \in C^{-\infty}(\tilde L_0)$. By \eqref{eq_invariance_rho} and the invariance of $\phi_L$ we obtain that  
\begin{displaymath}
 g_\alpha^*f=\psi_\alpha(\xi)^{-\frac{n-k}{2}}f,
\end{displaymath}
and applying Proposition \ref{prop:uniqueness_of_homogeneous_extension} shows that outside of $\xi_0$, $f$ must coincide with a multiple of $\sigma^{-\frac{n-k}{2}}$, that is $c_+=c_-$.

On the other hand, since $(n-k)$ is even, $\beta=-1$ if and only if $\nu=\frac{n-k}{2}-1$. By \eqref{eq_inequality_from_lemma} it follows that the space of $(G,\det)$-invariant generalized primitive forms of degree $(k,n-k-1)$ supported on the light cone is at most $1$-dimensional. 

In both cases we find $\dim \Curv_k^{-\infty,\OO(p,q)} = \dim \Omega_{p,-\infty}^{k,n-k-1}(\p_V)^{tr,G,\det}\leq 2$.
\endproof

\proof[Proof of Theorem \ref{mainthm_uniqueness_functors}]

Consider a Weyl functor $\Lambda:\mathbf{\Psi Met}\to  \mathbf {GCrv}$. Fix $p,q>0$. From Proposition \ref{prop_flat_case} we obtain that 
\begin{displaymath} 
\Lambda^{\R^{p,q}}=\sum_{k=0}^{p+q} a_k \Lambda_k^{\R^{p,q}}+b_k \bar \Lambda_k^{\R^{p,q}}
\end{displaymath}
for some constants $a_k,b_k$. Since the $\Lambda_k,\bar \Lambda_k$ are linearly independent and satisfy a Weyl principle, the $a_k,b_k$ are unique and independent of the choice of $\R^{p,q}$ provided that $k<p+q$. Then $\Lambda^{\R^{p',q'}}=\sum_{k=0}^\infty a_k \Lambda_k^{\R^{p',q'}}+b_k \bar \Lambda_k^{\R^{p',q'}}$ on each pseudo-Euclidean space $\R^{p',q'}$. By functoriality and the pseudo-Riemannian Nash embedding theorem \cite{clarke70}, we then have on each pseudo-Riemann manifold $M$
\begin{displaymath} 
\Lambda^{M}=\sum_{k=0}^\infty a_k \Lambda_k^M+b_k \bar \Lambda_k^M.
\end{displaymath} 
\endproof


\section{A K\"unneth-type formula for Lipschitz-Killing curvature measures}

\subsection{Disintegration of curvature measures}

We start with a general proposition.

\begin{Proposition} \label{prop_disintegration_smooth}
Let $X_1,X_2$ be smooth manifolds and $X:=X_1 \times X_2$.
Let $\Psi \in \mathcal{C}^\infty(X), \phi_2 \in \mathcal{V}_c^\infty(X_2), f_2 \in C_c^\infty(X_2)$. Then  there exists a unique smooth curvature measure $\tilde \Psi \in \mathcal{C}^\infty(X_1)$ such that 
\begin{equation} \label{eq_disintegration_curvature_measure_smooth}
\tilde \Psi(\phi_1,f_1)=\Psi(\phi_1 \boxtimes \phi_2,f_1 \boxtimes f_2), \quad \phi_1 \in \mathcal{V}_c^\infty(X_1), f_1 \in C^\infty(X_1).
\end{equation}
\end{Proposition}

\proof
Uniqueness is clear. To prove existence, we use the notations and maps from \cite[Section 4.1]{alesker_bernig_convolution}. The relevant diagram is

\begin{displaymath}
\xymatrix{\p_{X_1}
& \p_{X_1} \times \p_{X_2} \ar@{->>}[l]_-{q_1} \ar@{->>}[r]^-{q_2} & \p_{X_2}\\
& \hat \p_X \ar@{->>}[u]^\Phi \ar@{->>}[d]_F & \\
\p_{X_1} \times X_2 \ar@{->>}[d]_{p_1} \ar@{^{(}->}[r]^-{i_1} & \p_X \ar@{->>}[dd]_{\pi_X} &  X_1 \times \p_{X_2}
\ar@{_{(}->}[l]_-{i_2} \ar@{->>}[d]_{p_2}\\
\p_{X_1} \ar@{->>}[d]_{\pi_{X_1}} & & \p_{X_2} \ar@{->>}[d]_{\pi_{X_2}}\\
X_1 & X=X_1 \times X_2 \ar@{->>}[l]_-{\tilde p_1} \ar@{->>}[r]^-{\tilde p_2} & X_2}
\end{displaymath}

Let $\phi_i \in \mathcal{V}_c^\infty(X_i), i=1,2$ be represented by forms $\phi_i \in \Omega^{n_i}(X_i), \omega_i \in \Omega^{n_i-1}(\p_{X_i})$. According to \cite{alesker_intgeo, alesker_bernig_convolution}, the exterior product $\phi_1 \boxtimes \phi_2 \in \mathcal{V}^{-\infty}(X_1 \times X_2)$ is represented by generalized forms $\phi \in \Omega^{n_1+n_2}_{-\infty}(X_1 \times X_2),\omega \in \Omega^{n_1+n_2-1}_{-\infty}(\p_X)$ such that 
\begin{align}
 (\pi_{X})_*\omega & =\tilde p_1^*(\pi_{X_1})_*\omega_1 \cdot \tilde p_2^*(\pi_{X_2})_*\omega_2 \label{eq_ext_product_first_term} \\
 a^*(D\omega+\pi_X^*\phi) & = F_*\Phi^*\left[q_1^*a^*(D\omega_1+\pi_{X_1}^*\phi_1) \wedge q_2^*a^*(D\omega_2+\pi_{X_2}^*\phi_2)\right] \nonumber \\
 & \quad + (\tilde p_1 \circ \pi_X)^* (\pi_{X_1})_*\omega_1 \wedge i_{2*}p_2^*a^*(D\omega_2+\pi_{X_2}^*\phi_2) \nonumber \\
 & \quad +
i_{1*}p_1^*a^*(D\omega_1+\pi_{X_1}^*\phi_1) \wedge (\tilde p_2 \circ \pi_X)^* (\pi_{X_2})_*\omega_2. \nonumber
\end{align}

The second equation implies that
\begin{align} 
D\omega+\pi_X^*\phi & = F_*\Phi^*\left[q_1^*(D\omega_1+\pi_{X_1}^*\phi_1) \wedge q_2^*(D\omega_2+\pi_{X_2}^*\phi_2)\right] \nonumber \\
 & \quad + (-1)^{n_1} (\tilde p_1 \circ \pi_X)^* (\pi_{X_1})_*\omega_1 \wedge i_{2*}p_2^*(D\omega_2+\pi_{X_2}^*\phi_2) \nonumber \\
 & \quad +
(-1)^{n_2}i_{1*}p_1^*(D\omega_1+\pi_{X_1}^*\phi_1) \wedge (\tilde p_2 \circ \pi_X)^* (\pi_{X_2})_*\omega_2. \label{eq_ext_product_second_term}
\end{align}
The signs come from the fact that the degree of the antipodal map is given by $(-1)^{n_1}$ on $\p_{X_1} \times X_2$; by  $(-1)^{n_2}$ on $X_1 \times \p_{X_2}$; and $(-1)^{n_1+n_2}$ on $\p_X$ and on $\hat \p_X$.

Let $\Psi \in \mathcal{C}^\infty(X_1 \times X_2)$ be a smooth curvature measure, given by forms $\rho \in \Omega^{n_1+n_2}(X_1 \times X_2),\eta \in \Omega^{n_1+n_2-1}(\p_X)$. Then 
\begin{displaymath}
 \Psi(\phi_1 \boxtimes \phi_2,f_1 \boxtimes f_2)=\int_X \tilde p_1^*f_1 \cdot \tilde p_2^*f_2 \cdot (\pi_{X})_*\omega \wedge \rho + \int_{\p_X}  \pi_X^* \tilde p_1^*f_1 \cdot \pi_X^* \tilde p_2^*f_2 \cdot (D\omega+\pi_X^*\phi) \wedge \eta.
\end{displaymath}

Using \eqref{eq_ext_product_first_term}, we see that the first summand is given by 
\begin{displaymath}
 \int_X \tilde p_1^*(f_1 \cdot (\pi_{X_1})_*\omega_1) \cdot \tilde p_2^*(f_2 \cdot (\pi_{X_2})_*\omega_2) \wedge \rho=\int_{X_1} f_1 \cdot (\pi_{X_1})_*\omega_1 \wedge (\tilde p_1)_*\left[\tilde p_2^*(f_2  \cdot (\pi_{X_2})_*\omega_2) \wedge \rho\right].
\end{displaymath}

According to \eqref{eq_ext_product_second_term}, the second summand splits as the sum $T_1+T_2+T_3$ where 
\begin{align*}
 T_1 & = \int_{\p_X}  F_*\Phi^*\left[q_1^*(D\omega_1+\pi_{X_1}^*\phi_1) \wedge q_2^*(D\omega_2+\pi_{X_2}^*\phi_2)\right] \wedge \pi_X^* \tilde p_1^*f_1 \cdot \pi_X^* \tilde p_2^*f_2 \wedge \eta\\
 & = \int_{\p_{X_1} \times \p_{X_2}}  q_1^*(D\omega_1+\pi_{X_1}^*\phi_1) \wedge q_2^*(D\omega_2+\pi_{X_2}^*\phi_2) \wedge q_1^*\pi_{X_1}^*f_1 \wedge q_2^*\pi_{X_2}^*f_2 \wedge \Phi_*F^* \eta\\
  & = \int_{\p_{X_1}}  \pi_{X_1}^*f_1 \cdot (D\omega_1+\pi_{X_1}^*\phi_1)  \wedge (q_1)_*\left[q_2^*(D\omega_2+\pi_{X_2}^*\phi_2) \wedge q_2^*\pi_{X_2}^*f_2 \wedge \Phi_*F^* \eta\right],\\
T_2 & = (-1)^{n_1} \int_{\p_X}  \pi_X^* \tilde p_1^*f_1 \cdot \pi_X^* \tilde p_2^*f_2 \cdot (\tilde p_1 \circ \pi_X)^* (\pi_{X_1})_*\omega_1 \wedge i_{2*}p_2^*(D\omega_2+\pi_{X_2}^*\phi_2) \wedge \eta\\
& = (-1)^{n_1} \int_{X_1}  f_1 \cdot (\pi_{X_1})_*\omega_1 \wedge (\tilde p_1 \circ \pi_X)_*\left[\pi_X^* \tilde p_2^*f_2 \cdot i_{2*}p_2^*(D\omega_2+\pi_{X_2}^*\phi_2) \wedge \eta\right],\\
T_3 & = (-1)^{n_2} \int_{\p_X}  \pi_X^* \tilde p_1^*f_1 \cdot \pi_X^* \tilde p_2^*f_2 \cdot i_{1*}p_1^*(D\omega_1+\pi_{X_1}^*\phi_1) \wedge (\tilde p_2 \circ \pi_X)^* (\pi_{X_2})_*\omega_2 \wedge \eta\\
& = (-1)^{n_2} \int_{\p_{X_1} \times X_2}  i_1^*\left[\pi_X^* \tilde p_1^*f_1 \cdot \pi_X^* \tilde p_2^*f_2 \cdot (\tilde p_2 \circ \pi_X)^* (\pi_{X_2})_*\omega_2 \wedge \eta\right] \wedge p_1^*(D\omega_1+\pi_{X_1}^*\phi_1)\\
& = (-1)^{n_2+n_1n_2} \int_{\p_{X_1} \times X_2}  p_1^*(D\omega_1+\pi_{X_1}^*\phi_1) \wedge i_1^*\left[\pi_X^* \tilde p_1^*f_1 \cdot \pi_X^* \tilde p_2^*f_2 \cdot (\tilde p_2 \circ \pi_X)^* (\pi_{X_2})_*\omega_2 \wedge \eta\right] \\
& = (-1)^{n_2+n_1n_2} \int_{\p_{X_1}}  \pi_{X_1}^*f_1 \cdot (D\omega_1+\pi_{X_1}^*\phi_1) \wedge (p_1)_*i_1^*\left[\pi_X^* \tilde p_2^*f_2 \cdot (\tilde p_2 \circ \pi_X)^* (\pi_{X_2})_*\omega_2 \wedge \eta\right].
\end{align*}

We set 
\begin{align}
\tilde \phi & := (\tilde p_1)_*\left[\tilde p_2^*f_2  \cdot \tilde p_2^*(\pi_{X_2})_*\omega_2 \wedge \rho\right]   \nonumber \\
& \quad +(-1)^{n_1}(\tilde p_1 \circ \pi_X)_*\left[\pi_X^* \tilde p_2^*f_2 \cdot i_{2*}p_2^*(D\omega_2+\pi_{X_2}^*\phi_2) \wedge \eta\right] \in \Omega^{n_1}(X_1) \label{eq_def_tilde_phi}\\
 \tilde \omega & := (q_1)_*\left[q_2^*(D\omega_2+\pi_{X_2}^*\phi_2) \wedge q_2^*\pi_{X_2}^*f_2 \wedge \Phi_*F^* \eta\right] \nonumber \\
 & \quad + (-1)^{n_1+n_1n_2} (p_1)_*i_1^*\left[\pi_X^* \tilde p_2^*f_2 \cdot (\tilde p_2 \circ \pi_X)^* (\pi_{X_2})_*\omega_2 \wedge \eta\right] \in \Omega^{n_1-1}(\p_{X_1}). \label{eq_def_tilde_omega}
\end{align}
Since $p_1,\tilde p_1,q_1$ are submersions, it follows that $\tilde \omega$ and the first summand of $\tilde \phi$ are smooth. For the second summand, this is not immediate because of the push-forward under the map $i_2$ which is not a submersion. But the restriction of $\tilde p_1 \circ \pi_X$ to the image of $i_2$ is obviously a submersion, hence this term is smooth as well.

We thus see that the curvature measure $\tilde \Psi:=[\tilde \phi,\tilde \omega] \in \mathcal{C}^\infty(X_1)$ satisfies \eqref{eq_disintegration_curvature_measure_smooth}. 
\endproof

We need a version of the previous proposition for generalized curvature measures. In this case, the left hand side of \eqref{eq_disintegration_curvature_measure_smooth} is well-defined only if $\Psi$ satisfies an additional condition of transversality, compare also \cite[Section 4]{alesker_bernig}.

Recall from \cite{alesker_bernig_convolution} the set $\mathcal M=\mathcal M_1 \cup \mathcal M_2 \subset \p_X$ with
\begin{align*}
\mathcal{M}_1 & := \mathrm {im} i_1=\{(x_1,x_2,[\xi_1:0]), x_1 \in X_1,x_2 \in X_2,\xi_1 \in T^*_{x_1}X_1 \setminus 0\},\\
\mathcal{M}_2 & := \mathrm {im} i_2=\{(x_1,x_2,[0:\xi_2]), x_1 \in X_1,x_2 \in X_2,\xi_2 \in T^*_{x_2}X_2 \setminus 0\}.
\end{align*}

\begin{Definition} \label{def_transversal_curvature_measure}
Let $X_1,X_2$ be smooth manifolds. Then a generalized curvature measure $\Psi \in \mathcal{C}^{-\infty}(X_1 \times X_2)$ is called \emph{transversal}, if it belongs to 
\begin{displaymath}
\bigcup_\Gamma \mathcal{C}^{-\infty}_{T^*X \setminus \underline{0},\Gamma}(X_1 \times X_2),
\end{displaymath}
where $\Gamma$ runs over all closed subsets in $T^*\p_X \setminus \underline{0}$ that are disjoint from $T^*_{\mathcal M}\p_X$.
\end{Definition}

By \cite[Proposition 4.1]{alesker_bernig_convolution}, a transversal curvature measure can be applied to a pair $(\phi_1 \boxtimes \phi_2,f_1 \boxtimes f_2)$ with $\phi_i \in \mathcal{V}^\infty_c(X_i), f_i \in C_c^\infty(X_i)$.

\begin{Proposition} \label{prop_disintegration_generalized}
Let $X_1,X_2$ be smooth manifolds. Let $\Psi \in \mathcal{C}^{-\infty}(X_1 \times X_2)$ be transversal and $\phi_2 \in \mathcal{V}_c^\infty(X_2), f_2 \in C_c^\infty(X_2)$. Then  there exists a unique generalized curvature measure $\tilde \Psi \in \mathcal{C}^{-\infty}(X_1)$ such that 
\begin{equation} \label{eq_disintegration_curvature_measure_generalized}
\tilde \Psi(\phi_1,f_1)=\Psi(\phi_1 \boxtimes \phi_2,f_1 \boxtimes f_2), \quad \phi_1 \in \mathcal{V}_c^\infty(X_1), f_1 \in C_c^\infty(X_1).
\end{equation}
\end{Proposition}

\proof
Let $\Psi \in \mathcal{C}^{-\infty}_{T^*X \setminus \underline{0},\Gamma}(X_1 \times X_2)$, where $\Gamma$ is disjoint from $T^*_{\mathcal M}\p_X$. Let $\Psi$ be represented by generalized forms $\rho \in \Omega^{n_1+n_2}_{-\infty}(X_1 \times X_2), \eta \in \Omega^{n_1+n_2-1}_\Gamma(\p_X)$.

As in the previous proof, we are going to define the forms $\tilde \phi \in \Omega^{n_1}_{-\infty}(X), \tilde \omega \in \Omega^{n_1-1}_{-\infty}(\p_X)$ by \eqref{eq_def_tilde_phi} and \eqref{eq_def_tilde_omega}. 

We have to check that this is possible. For $\tilde \phi$ we note that the wave front set of ${i_2}_{*}p_2^*(D\omega_2+\pi_{X_2}^*\phi_2)$ is contained in $T^*_{\mathcal M_2}\p_X \subset T^*_{\mathcal M}\p_X$, hence the wedge product with $\eta$ is defined, and then $\tilde \phi$ is well-defined.

The first term in the definition of $\tilde \omega$ is well-defined since $F$ is a submersion. The second term is well-defined since $i_1$ is transversal to $\Gamma$. 

From the same arguments as in the previous proof, we see that the generalized curvature measure $\tilde \Psi:=[\tilde \phi,\tilde \omega]$ satisfies \eqref{eq_disintegration_curvature_measure_generalized}.
\endproof

\subsection{Proof of the K\"unneth-type formula}
\begin{Proposition} \label{prop_product_lc}
Let $(M_i,Q_i)$ be pseudo-Riemannian manifolds and $\psi_i \in \mathcal{V}^{-\infty}(M_i)$ for $i=1,2$ be LC-transversal. Then $\psi_1 \boxtimes \psi_2 \in \mathcal{V}^{-\infty}(M_1 \times M_2)$ is LC-transversal. More precisely, if $\WF(\psi_i) \subset (\Lambda_i,\Gamma_i)$ and $\Gamma_i \cap N^*(\LC^*_{M_i})=\emptyset$ for $i=1,2$, then there exists $(\Lambda,\Gamma)$ with $\Gamma \cap N^*(\LC^*_{M_1 \times M_2}) = \emptyset$ and such that 
\begin{equation}
\boxtimes: \mathcal{V}^{-\infty}_{\Lambda_1,\Gamma_1}(M_1) \times \mathcal{V}^{-\infty}_{\Lambda_2,\Gamma_2}(M_2) \to \mathcal{V}^{-\infty}_{\Lambda,\Gamma}(M_1 \times M_2)
\end{equation}
is continuous.    
\end{Proposition}

\proof
We prove the statement in the case of $LC$-transversal submanifolds $X_1,X_2$. 

Let $Q_1,Q_2$ be the metrics on $M_1,M_2$ and $Q:=Q_1 \oplus Q_2$. Since they are non-degenerate, the manifolds 
\begin{align*}
 \LC^*_{M_i} & = \{(x,[\xi]) \in \p_{M_i}: Q_i(\xi)=0\}, \quad i=1,2,\\
 \LC^*_{M_1 \times M_2} & = \{(x_1,x_2,[\xi_1:\xi_2]) \in \p_{M_1 \times M_2}: Q_1(\xi_1)+Q_2(\xi_2)=0\}
\end{align*}
are of codimension $1$. 

The conormal bundle of $X_1 \times X_2$ is given by 
\begin{displaymath}
 N^*(X_1 \times X_2) =\{(x_1,x_2,[\xi_1:\xi_2]) \in \p_{M_1 \times M_2}: \xi_i=0 \text{ or } (x_i,[\xi_i]) \in N^*{X_i}, \forall i=1,2\}.
\end{displaymath}

Let $(x_1,x_2,[\xi_1:\xi_2]) \in N^*(X_1 \times X_2) \cap \LC^*_{M_1 \times M_2}$. We consider two cases. If $Q_1(\xi_1) \neq 0$, then $\xi_1 \neq 0$ and the curve $c(t):=(x_1,x_2,[t\xi_1:\xi_2])$ stays inside $N^*(X_1 \times X_2)$. But the vector $c'(1)$ is not tangent to $\LC^*_{M_1 \times M_2}$.  

If $Q_1(\xi_1)=0$, then also $Q_2(\xi_2)=0$. Suppose that $\xi_1 \neq 0$. Since $X_1$ is LC-regular, there exists some tangent vector $w_1 \in T_{(x_1,[\xi_1])} N^*X_1$ with $w_1 \not\in T_{(x_1,[\xi_1])} \LC^*_{M_1}$. Then $(w_1,0)$ is tangent to $N^*(X_1 \times X_2)$, but not tangent to $\LC^*_{M_1 \times M_2}$. Finally, if $\xi_1=0$ then $\xi_2 \neq 0$ and we may argue as before, using that $X_2$ is LC-regular. 

We omit the proof of the general case of the statement. In this case, a careful analysis of the wave front set of the exterior product $\psi_1 \boxtimes \psi_2$ is needed. By \cite[Proposition 4.1]{alesker_bernig_convolution} it consists of the union of three sets, and each of these sets can be shown to be disjoint from the wave front set of $\LC^*_{M_1 \times M_2}$ by arguments which are similar to the given ones.
\endproof

\proof[Proof of Theorem \ref{mainthm_kuenneth}]

By \cite[Proposition 7.1]{bernig_faifman_solanes}, we have 
\begin{displaymath}
\Lambda_k^{M_1 \times M_2} \in \mathcal{C}^{-\infty}_{\emptyset, N^*(\LC^*_{M_1 \times M_2})}(M_1 \times M_2).
\end{displaymath} 
We claim that 
\begin{displaymath}
N^*(\LC^*_{M_1 \times M_2}) \cap T^*_{\mathcal M}\p_{M_1 \times M_2}=\emptyset,
\end{displaymath}
which implies that $\Lambda_k^{M_1 \times M_2}$ is transversal in the sense of Definition \ref{def_transversal_curvature_measure}. To prove the claim, take an element in $\LC^*_{M_1 \times M_2} \cap \mathcal M$. Without loss of generality, suppose that it belongs to $\mathcal M_1$, i.e. it is of the form $(x_1,x_2,[\xi_1:0])$ with $x_1 \in M_1,x_2 \in M_2$ and $\xi_1$ a null-covector of $M_1$. The claim now follows from the fact that $0$ is a regular value of $Q_1 \in C^\infty(TM_1 \setminus \underline{0})$. 
 
Since $\Lambda_k^{M_1 \times M_2}$ is a generalized curvature measure on $M_1 \times M_2$, we may integrate a smooth function of the form $f_1 \boxtimes f_2 \in C_c^\infty(M_1 \times M_2)$ with $f_i \in C_c^\infty(M_i),i=1,2$ and the result is a generalized valuation on $M_1 \times M_2$.

By Propositions \ref{prop_product_lc} and \cite[Proposition 7.1]{bernig_faifman_solanes}, it can be paired with an exterior product $\phi_1 \boxtimes \phi_2$ of smooth valuations $\phi_i \in \mathcal{V}^\infty_c(M_i), i=1,2$. 

By Proposition \ref{prop_disintegration_generalized},
\begin{displaymath}
 (\phi_1,f_1) \mapsto \Lambda_k^{M_1 \times M_2}(\phi_1 \boxtimes \phi_2,f_1 \boxtimes f_2)
\end{displaymath}
is a generalized curvature measure on $M_1$, which will be denoted by $\Lambda_k^{M_1 \times M_2}(\bullet \boxtimes \phi_2,\bullet \boxtimes f_2)$.

If $i:M_1 \to \tilde M_1$ is an isometric embedding, then $i \times \mathrm{id}: M_1 \times M_2 \to \tilde M_1 \times M_2$ is also an isometric embedding. By Weyl's principle \cite[Theorem D]{bernig_faifman_solanes} we have 
\begin{displaymath}
(i \times \mathrm{id})^* \Lambda_k^{\tilde M_1 \times M_2}(\bullet \boxtimes \phi_2, \bullet \boxtimes f_2)=\Lambda_k^{M_1 \times M_2}(\bullet \boxtimes \phi_2,\bullet \boxtimes f_2). 
\end{displaymath}

Hence $ M_1 \mapsto \Lambda_k^{M_1 \times M_2}(\bullet \boxtimes \phi_2,\bullet \boxtimes f_2)$ is a Weyl functor. By Theorem \ref{mainthm_uniqueness_functors}, there are constants $a_{k_1}^{M_2,k}(\phi_2,f_2)$ and $b_{k_1}^{M_2,k}(\phi_2,f_2), k=0,1,\ldots$ which do not depend on $M_1$ such that 
\begin{equation} \label{eq_product1}
\Lambda_k^{M_1 \times M_2}(\bullet \boxtimes \phi_2, \bullet \boxtimes f_2\rangle=\sum_{k_1=0}^\infty a^{M_2,k}_{k_1}(\phi_2,f_2) \Lambda_{k_1}^{M_1}+b^{M_2,k}_{k_1}(\phi_2,f_2) \overline{\Lambda}_{k_1}^{M_1}
\end{equation}
for all $M_1$.

By induction on $k_1$ we will show that for each fixed $M_2$, 
\begin{displaymath}
 (\phi_2,f_2) \mapsto a_{k_1}^{M_2,k}(\phi_2,f_2),b_{k_1}^{M_2,k}(\phi_2,f_2)
\end{displaymath}
are Weyl functors. Suppose that this is true for all $k_1<n$. 

Recall that for an $n$-dimensional pseudo-Riemannian manifold $M_1$ of signature $(p,q)$ we have that $\Lambda_{n}^{M_1}=\mathbf i^q \vol_{M_1}$. We can therefore choose two $n$-dimensional pseudo-Riemannian manifolds $M_1,M_1'$ and $f_1 \in C_c^\infty(M_1), f_1' \in C_c^\infty(M_1'), \phi_1 \in \mathcal{V}_c^\infty(M_1),\phi_1' \in \mathcal{V}_c^\infty(M_1')$ such that the matrix 
\begin{displaymath}
\left(\begin{matrix}
       \Lambda_n^{M_1}(\phi_1,f_1) & \bar \Lambda_n^{M_1}(\phi_1,f_1)\\
        \Lambda_n^{M_1'}(\phi_1',f_1') & \bar \Lambda_n^{M_1'}(\phi_1',f_1')
      \end{matrix}
 \right)
\end{displaymath}
is non singular. 

By \eqref{eq_product1} we find that 
\begin{align*}
 \Lambda_{n}^{M_1}&(\phi_1,f_1) a_{n}^{M_2,k} +\bar \Lambda_n^{M_1}(\phi_1,f_1) b_n^{M_2,k}\\
  & =\Lambda_k^{M_1 \times M_2}(\phi_1 \boxtimes \bullet,f_1 \boxtimes \bullet)-\sum_{k_1=0}^{n-1} \left(\Lambda_{k_1}^{M_1}(\phi_1,f_1) a_{k_1}^{M_2,k}+\bar \Lambda_{k_1}^{M_1}(\phi_1,f_1) b_{k_1}^{M_2,k}\right)\\
 \Lambda_{n}^{M_1'}&(\phi_1',f_1') a_{n}^{M_2,k} +\bar \Lambda_n^{M_1'}(\phi_1',f_1') b_n^{M_2,k} \\
 & =\Lambda_k^{M_1' \times M_2}(\phi_1' \boxtimes \bullet,f_1' \boxtimes \bullet)-\sum_{k_1=0}^{n-1} \left( \Lambda_{k_1}^{M_1'}(\phi_1',f_1') a_{k_1}^{M_2,k}+\bar \Lambda_{k_1}^{M_1'}(\phi_1',f_1') b_{k_1}^{M_2,k}\right).
\end{align*}

 The first term on the right hand side of each equation is a Weyl functor on $M_2$ by an argument as above. The second term on the right hand side is a Weyl functor by the induction hypothesis. Solving this system of equations yields that $M_2 \mapsto a_n^{M_2,k}, M_2 \mapsto b_n^{M_2,k}$ are Weyl functors as claimed.    

We thus may write 
\begin{align*}
 a_{k_1}^{M_2,k} & =\sum_{k_2=0}^\infty a_{k_1,k_2}^k \Lambda_{k_2}^{M_2}+\hat a_{k_1,k_2}^k  \bar \Lambda_{k_2}^{M_2},\\
 b_{k_1}^{M_2,k} & =\sum_{k_2=0}^\infty b_{k_1,k_2}^k \Lambda_{k_2}^{M_2}+\hat b_{k_1,k_2}^k  \bar \Lambda_{k_2}^{M_2},
\end{align*}
with scalars $a_{k_1,k_2}^k,\hat a_{k_1,k_2}^k, b_{k_1,k_2}^k,\hat b_{k_1,k_2}^k$ which do not depend on $M_1$ and $M_2$. Then 
\begin{align*} 
 \Lambda_k^{M_1 \times M_2}&(\phi_1 \boxtimes \phi_2,\bullet   \boxtimes \bullet) \\
 & =\sum_{k_1,k_2} a_{k_1,k_2}^k \Lambda_{k_1}^{M_1}(\phi_1,\bullet) \boxtimes \Lambda_{k_2}^{M_2}(\phi_2,\bullet)+\hat a_{k_1,k_2}^k \Lambda_{k_1}^{M_1}(\phi_1,\bullet) \boxtimes \bar \Lambda_{k_2}^{M_2}(\phi_2,\bullet) \nonumber\\
 & \quad +b_{k_1,k_2}^k \bar \Lambda_{k_1}^{M_1}(\phi_1,\bullet) \boxtimes \Lambda_{k_2}^{M_2}(\phi_2,\bullet)+\hat b_{k_1,k_2}^k \bar \Lambda_{k_1}^{M_1}(\phi_1,\bullet) \boxtimes \bar \Lambda_{k_2}^{M_2}(\phi_2,\bullet). 
\end{align*}

Using the scaling property of the Lipschitz-Killing curvature measures \cite[Proposition 7.10]{bernig_faifman_solanes} we find that $a_{k_1,k_2}^k=\hat a_{k_1,k_2}^k=b_{k_1,k_2}^k=\hat b_{k_1,k_2}^k=0$ unless $k_1+k_2=k$. 

It remains to determine $a_{k_1,k_2}^{k_1+k_2},\hat a_{k_1,k_2}^{k_1+k_2},b_{k_1,k_2}^{k_1+k_2},\hat b_{k_1,k_2}^{k_1+k_2}$. Take pseudo-Riemannian manifolds $M_1,M_2$ of dimensions $k_1$ and $k_2$ respectively. Then, whatever the signatures are, we have $\vol_{M_1 \times M_2}=\vol_{M_1} \boxtimes \vol_{M_2}$. Taking riemannian and lorentzian signature metrics on $M_1$ and $M_2$ and using \cite[Proposition 7.7.]{bernig_faifman_solanes}, we obtain a system of four linear equations with the unique solution $a_{k_1,k_2}^{k_1+k_2}=1,\hat a_{k_1,k_2}^{k_1+k_2}=b_{k_1,k_2}^{k_1+k_2}=\hat b_{k_1,k_2}^{k_1+k_2}=0$.

We thus have
\begin{equation}
\Lambda_k^{M_1 \times M_2}(\phi_1 \boxtimes \phi_2,\bullet  \times \bullet) =\sum_{k_1+k_2=k} \Lambda_{k_1}^{M_1}(\phi_1,\bullet) \boxtimes \Lambda_{k_2}^{M_2}(\phi_2,\bullet)
\label{eq_product_property} 
\end{equation}
for all $\phi_1 \in \mathcal{V}_c^\infty(M_1), \phi_2 \in \mathcal{V}_c^\infty(M_2)$.
 
Finally, let $A_i  \in \mathcal P(M_i),i=1,2$ be LC-transversal. By Lemma \ref{lem:lc_transversal_implies_lc_transversal}, $\chi_{A_i}$ are LC-transversal generalized valuations. Then by Proposition \ref{prop_product_lc} it follows that $\chi_{A_1 \times A_2}=\chi_{A_1}\boxtimes \chi_{A_2}$ is LC-transversal. Set $(\Lambda_i,\Gamma_i):=\WF(\chi_{A_i})$. By assumption, we have
\begin{displaymath}
 \Gamma_i \cap N^* (\LC^*_{M_i} )= \emptyset, i=1,2. 
\end{displaymath}

By Proposition \ref{prop_product_lc} there exist some $(\Lambda,\Gamma)$ with 
\begin{displaymath}
 \Gamma \cap N^*( \LC^*_{M_1 \times M_2})=\emptyset,
\end{displaymath}
such that the exterior product $\boxtimes: \mathcal{V}^{-\infty}_{\Lambda_1,\Gamma_1}(M_1) \times \mathcal{V}^{-\infty}_{ \Lambda_2,\Gamma_2}(M_2) \to \mathcal{V}^{-\infty}_{\Lambda,\Gamma}(M_1 \times M_2)$ is continuous. 

Take a sequence $\phi_i^m \in \mathcal{V}_c^\infty(M_i)$ which converges to $\chi_{A_i}$ in $\mathcal{V}_{\Lambda_i,\Gamma_i}^{-\infty}(M_i)$. The sequence $\phi_1^m \boxtimes \phi_2^m$ then converges in $\mathcal{V}_{\Lambda,\Gamma}^{-\infty}(M_1 \times M_2)$ to $\chi_{A_1} \boxtimes \chi_{A_2}$.

By \cite[Proposition 7.1]{bernig_faifman_solanes} we may take the limit in Equation \eqref{eq_product_property} which yields the statement. 
\endproof

\section{Generalized valuations and curvature measures on isotropic space forms}

The aim of this section is to prove Theorem \ref{mainthm:classification_spq}. The following special case of valuations on pseudo-Euclidean spaces $\R^{p,q}$ was considered in \cite{alesker_faifman} for $\min(p,q)=1$ and in \cite{bernig_faifman_opq} in general.

\begin{Theorem}\label{thm:rpq_classification}
The space of translation and $\OO(p,q)$-invariant generalized valuations on $\R^{p,q}$ coincides with $\mathcal {LK}(\R^{p,q})$.
\end{Theorem}

The corresponding statement for curvature measures is Proposition \ref{prop_flat_case}. To extend this statement to isotropic space forms, we will need the following technical lemma.

\begin{Lemma}\label{lemma:dual_sections}
Let $M$ be a smooth manifold, $G$ a locally compact Lie group acting on $M$ transitively, and $\mathcal E$ a $G$-equivariant linear bundle over $M$, such that the fiber $\mathcal E|_x$ is a smooth Fr\'echet module of $\Stab(x)$. Then the $G$-invariants of the dual space $\Gamma_c^{\infty}(M, \mathcal E)^\ast$ coincide with $\Gamma^{\infty}(M, \mathcal E^*\otimes|\omega_M|)^G$.
\end{Lemma}

\proof 
Take $s \in \left(\Gamma_c^{\infty}(M, \mathcal E)^\ast\right)^G$. For $\phi\in\Gamma_c^\infty(M, \mathcal E)$ one can define $s\cdot \phi \in\mathcal M^{-\infty}(M)$ by setting $\int f\cdot d(s\cdot \phi):=s(f\phi)$ for all $f \in C^\infty_c(M)$.
The theorem of Dixmier-Malliavin \cite{dixmier_malliavin} (see \cite{casselman_notes_dixmiermalliavin} for a simpler proof) implies that $\phi$ can be represented as $\phi=\sum_{j=1}^N \int_G w_j(g)g^*\phi_jdg$ with $w_j \in C^\infty_c(G), \phi_j \in \Gamma^\infty(M,\mathcal E)$. Then by $G$-invariance of $s$,
\begin{displaymath}
s\cdot \phi =  \sum_{j=1}^N \int_G w_j(g)s\cdot g^*\phi_jdg= \sum_{j=1}^N \int_G w_j(g) g^*(s\cdot \phi_j)dg.
\end{displaymath}

It follows that $s\cdot\phi$ is a smooth measure. In particular, it defines a density on every tangent plane, $s\cdot\phi (x)\in\Dens(T_xM)$.
We next claim that $s\cdot\phi(x)$ only depends on $s$ and $\phi(x) \in \mathcal E|_x$. Indeed, if $\phi(x)=0$ we may represent $\phi=f\cdot\psi$ with  $\psi \in \Gamma_c^\infty(M,\mathcal E)$ and $f \in C_c^\infty(M), f(0)=0$. Then $(s\cdot \phi)=f\cdot (s\cdot \psi)$, and thus $(s\cdot \phi)(x)=0$. That is, $s\in \mathcal E^*|_x\otimes \Dens(T_xM)$, concluding the proof. 
\endproof

Consider the cosphere bundle $\pi:\p_M\to M$. The space $\Omega^{n-1}(\p_M)$ of differential forms has a natural filtration 
\begin{displaymath}
\Omega^{n-1}(\p_M)=\Omega_0^{n-1}(\p_M) \supset \Omega_1^{n-1}(\p_M) \supset \cdots \supset \Omega_{n-1}^{n-1}(\p_M) \supset \Omega_n^{n-1}(\p_M)=\{0\},
\end{displaymath}
where $\Omega^{n-1}_k(\p_M)=\langle \pi^*\Omega^k(M)\rangle_{n-1}$. Here $\langle A\rangle_{n-1}$ denotes the subspace of forms of degree $(n-1)$ in the ideal generated by $A$. This induces a filtration 
\begin{equation} \label{eq_filtration_smooth_curvature_measures}
\mathcal{C}^\infty(M)=\mathcal{C}^\infty_0(M) \supset \mathcal{C}^\infty_1(M) \supset \cdots \supset \mathcal{C}^\infty_{n}(M) \supset \mathcal{C}^\infty_{n+1}(M)=\{0\},
\end{equation}
where $\mathcal{C}^\infty_k(M)$ denotes the space of curvature measures $\Phi$ of the form $\Phi=[\omega,\phi]$ with $\omega \in \Omega^{n-1}_k(\p_M)$ for $0 \leq k \leq n-1$, and $\Phi=[0,\phi]$ with $\phi \in \Omega^n(M)$ for $k=n$. 

Denote by $\Curv_k^\infty(TM)$ the bundle over $M$ whose fiber at a point $x$ consists of the space $\Curv_k^\infty(T_xM)$. In \cite{solanes_wannerer}, a map (there denoted by $\Lambda_k'$) 
\begin{displaymath}
\sigma_k:\mathcal{C}^\infty_k(M) \to \Gamma^\infty(\Curv_k^\infty (TM))
\end{displaymath}
is constructed as follows. 

Fix $x \in M$ and take $\phi$ a local diffeomorphism $\phi:T_xM \to M$ with $\phi(0)=x, d\phi|_0=\mathrm{id}$. Then for $\Phi \in \mathcal{C}^\infty_k(M)$ 
\begin{displaymath}
\sigma_k \Phi(x):=\lim_{t \to 0} \frac{1}{t^k} (\phi \circ h_t)^*\Phi,
\end{displaymath}
where $h_t:T_xM \to T_xM$ is multiplication by $t$. 

In terms of forms, this map can be described as follows. If $k=n$, then $\Phi$ is a smooth measure on $M$, and evaluation at a point $x \in M$ gives us a Lebesgue measure on $T_xM$, hence an element of $\Curv_n^\infty(T_xM)$. 

If $k<n$, we may represent $\Phi=[\phi, \omega]$ with $\omega \in \Omega_k^{n-1}(\p_M)$.  Fix a point $(x,\xi) \in \p_M$. We use $d\pi^*$ to identify 
\begin{displaymath}
T_x^*M\simeq d\pi^*T_x^*M=T_\xi(\pi^{-1}(x))^\perp.  
\end{displaymath}

Consider the projection 
\begin{displaymath}
 \left\langle \largewedge^kT_x^*M \right\rangle_{n-1} \to \left\langle \largewedge^kT_x^*M \right\rangle_{n-1}/ \left\langle \largewedge^{k+1}T_x^*M\right\rangle_{n-1}= \largewedge^k T_x^*M \otimes \largewedge^{n-1-k}(T_{x,\xi}^*\p_M/T_x^*M),
\end{displaymath}
Noting that 
\begin{displaymath}
 T_{x,\xi}^*\p_M/T_x^*M\simeq T_\xi^*(\pi^{-1}(x)),
\end{displaymath}
this gives rise to a map
\begin{displaymath}
\tilde \sigma_k:\Omega^{n-1}_{k}(\mathbb P_M) \to \Gamma(\p_M, \largewedge^k T_x^*M \otimes \largewedge^{n-1-k}T_\xi^* \pi^{-1}(x)). 
\end{displaymath}
 
The fiber $\pi^{-1}(x)$ can be identified with the fiber of the cosphere bundle of $T_xM$, and so 
\begin{displaymath}
 \largewedge^k T_x^*M \otimes \largewedge^{n-1-k}T_\xi^* \pi^{-1}(x)=\largewedge^{k,n-1-k}T_{x,\xi}^*(\p_{T_xM}).
\end{displaymath}

Let $i\colon \pi^{-1}(x) \hookrightarrow \p_M$ be the inclusion. For $\omega\in \Omega_k^{n-1}(\p_M)$ we consider the restriction $i^*(\tilde\sigma_k\omega)$ (as a section of the pull-back bundle) and extend it to a translation-invariant form of $\p_{T_xM}$. We thus get a map
\begin{displaymath}
\bar \sigma_k:\Omega^{n-1}_{k}(\mathbb P_M)\to   \Gamma(M, \Omega^{k,n-1-k}(\mathbb P_{T_xM})^{\textrm{tr}}),
\end{displaymath}
which fulfills 
\begin{displaymath}
\sigma_k\Phi(x)=\sigma_k[\phi,\omega](x)=[0,\bar\sigma_k\omega(x)]. 
\end{displaymath}

Indeed, for $\omega=\pi^*\beta\wedge \eta$ with $\beta\in \Omega^k(M),\eta\in\Omega^{n-k-1}(\p_M)$ we have
\begin{align*}
 \tilde\sigma_k\omega(x,\xi)&=\beta_x \otimes r^*(\eta_\xi),
\end{align*}
where $r\colon T_\xi \pi^{-1}(x)\to T_\xi\p_M$ is the inclusion. On the other hand, if $\bar h_t,\bar\phi$ are the maps  on the cosphere bundles induced by $ h_t,\phi$, then for each $v\in T_xM$
\begin{align*}
 \lim_{t\to 0} t^{-k} (\bar h_t^*\bar\phi^*\omega)_{v,\xi}&=\lim_{t\to 0} t^{-k} (h_t^*\phi^*\beta)_v\wedge (\bar h_t^*\bar\phi^*\eta)_{v,\xi}=\beta_x\wedge r^*(\eta_\xi).
\end{align*}
Therefore, $\lim_{t\to 0} t^{-k} (\bar h_t^*\bar\phi^*\omega)$ is translation-invariant and equals $\bar\sigma_k\omega(x)$ as claimed.

The kernel of $\sigma_k$ is $\mathcal{C}^\infty_{k+1}(M)$, while $\sigma_k$ is evidently onto, hence we have an isomorphism 
\begin{displaymath}
\sigma_k: \mathcal{C}^\infty_k(M)/\mathcal{C}^\infty_{k+1}(M) \stackrel{\simeq}{\longrightarrow} \Gamma^\infty(\Curv_k^\infty(TM)).
\end{displaymath}

We now generalize these concepts to generalized curvature measures under an additional assumption. 

\begin{Definition}
 Let $M$ be a smooth manifold and $\Phi \in \mathcal{C}^{-\infty}(M)$. We define the space of \emph{horizontally smooth curvature measures} by 
\begin{displaymath}
\mathcal{C}^{hs}(M):=\bigcup_\Gamma \mathcal{C}_{\emptyset,\Gamma}^{-\infty}(M),
\end{displaymath}
where $\Gamma$ runs over all closed subsets in $T^*\p_M \setminus \underline{0}$ that are disjoint from $d\pi^*(T^*M) \setminus \underline{0}$. The topology is the inductive limit topology, in particular a sequence $\Phi_j \in \mathcal{C}^{hs}(M)$ converges to some $\Phi \in \mathcal{C}^{hs}(M)$ if there is some fixed $\Gamma$ such that $\Phi_j \in \mathcal{C}^{-\infty}_{\emptyset,\Gamma}(M)$ for all $j$ and  $\Phi_j \to \Phi$ in $\mathcal{C}^{-\infty}_{\emptyset,\Gamma}(M)$.
\end{Definition} 

Note that $\Omega_k^{n-1}(\p_M)$ equals the space of smooth sections of some finite-dimensional vector bundle over $\p_M$. We let $\Omega_k^{n-1,-\infty}(\p_M)$ be the corresponding space of generalized sections of the same vector bundle and define $\mathcal{C}_k^{-\infty}(M)$ as above. Then we have a filtration analogous to \eqref{eq_filtration_smooth_curvature_measures}.

\begin{Proposition} \label{prop_extension_solanes_wannerer_map}
 Let $M$ be a smooth manifold and $x \in M$. Then the map 
 \begin{displaymath}
  \sigma_k:\mathcal{C}_k^\infty(M) \to \Gamma^\infty(\Curv_k^\infty(TM)) 
 \end{displaymath}
 can be extended by continuity to a map 
 \begin{displaymath}
  \sigma_k: \mathcal{C}_k^{hs}(M) \to \Gamma(\Curv_k^{-\infty}(TM)), 
 \end{displaymath}
 whose kernel is $\mathcal{C}_{k+1}^{hs}(M)$. 
\end{Proposition}

\proof
Fix $x \in M$. A generalized form $\omega \in \Omega^{n-1,-\infty}_k(\p_M)$ can be restricted, as a section of the vector bundle with fiber $\largewedge^{n-1}T_{x,\xi}^*\p_M$, to the fiber $\pi^{-1}(x)$ under the assumption that $\WF(\omega)$ and $N^*\pi^{-1}(x)$ are disjoint, which is precisely what we get from the horizontal smoothness. Clearly the corresponding restriction map is continuous. Tracing the construction in the smooth case, we obtain a generalized, translation-invariant form on $\mathbb P_{T_xM}$ of bidegree $(k, n-1-k)$, and hence an element of $\Curv_k^{-\infty}(TM)$.  Using Proposition \ref{prop_kernel_generalized}, it is easy to check that the resulting map  vanishes on those forms which induce the trivial curvature measure, completing the proof.
\endproof

\begin{Proposition} \label{prop_invariant_implies_horizontally_smooth}
 Let $M$ be a smooth manifold, $G$ a Lie group that acts smoothly and transitively on $M$. Then 
 \begin{displaymath}
 \mathcal{C}^{-\infty}(M)^G \subset  \mathcal{C}^{hs}(M). 
 \end{displaymath}
\end{Proposition}

\proof
Let $X_1,\ldots,X_m$ be a basis of the Lie algebra $\mathfrak g$ and let $X_1^\#,\ldots,X^\#_m$ be the induced vector fields on $M$. Since $G$ acts transitively, we have that at each point $x \in M$, $X^\#_1|_x,\ldots,X^\#_m|_x$ generate the tangent space $T_xM$.  These vector fields induce in a natural way vector fields $\tilde X_1^\#,\ldots,\tilde X_m^\#$ on $\p_M$ such that $d\pi(\tilde X_i^\#)=X_i^\#$.

By Proposition \ref{prop_kernel_generalized}, a generalized curvature measure can be identified with the generalized sections of a certain finite-dimensional vector bundle $\mathcal E$ over $\p_M$. 

If $\Phi$ is $G$-invariant, then this section is $G$-invariant. Let $P_i$ denote the differential operator acting on sections of $\mathcal E$ corresponding to $\tilde X_i^\#$. By \cite[[Theorem 8.3.1]{hoermander_pde1} we have 
\begin{displaymath}
 \WF(\Phi) \subset \mathrm{Char}(P_i) \cup \underbrace{\WF(P_i\Phi)}_{=\emptyset}.
\end{displaymath}

The characteristic set of $P_i$ at a given point $(x,[\xi]) \in \p_M$ consists of all $\eta \in T^* \pi^{-1}(x) \setminus \{0\}$ that vanish on $\tilde X_i^\#$. On the other hand, the conormal bundle of $\pi^{-1}(x)$ consists of all such $\eta$ that vanish on all tangent vectors to $\pi^{-1}(x)$.

Now the vectors $\tilde X_1^\#,\ldots,\tilde X_m^\#$ and the tangent vectors to $\pi^{-1}(x)$ span $T_{(x,[\xi])}\p_M$, hence 
$\WF(\Phi) \cap N^* \pi^{-1}(x)=\emptyset$ as claimed. 
\endproof

\proof[Proof of Theorem \ref{mainthm:classification_spq}]
Let $M$ be an isotropic space form of signature $(p,q)$.  The isometry group $G:=\Isom(M)$ acts transitively on $M$. By the Weyl principle, the intrinsic volumes are invariant under isometries, i.e. $\mathcal{LK}(M) \subset \mathcal{V}^{-\infty}(M)^G$. We have to prove the opposite inclusion. 

Given $\mu \in \mathcal{V}^{-\infty}(M)^G$, take $k \in \{0,\ldots,n+1\}$ with $\mu \in \mathcal{W}_k^{-\infty}(M)$. If $k=n+1$, then $\mu=0$ and we are done. Suppose that $k<n+1$. 

We note that, by Proposition 7.3.2 of \cite{alesker_val_man4},
\begin{align*}
\mathcal W_k^{-\infty}(M)/\mathcal W_{k+1}^{-\infty}(M) & = \left(\mathcal W_{n-k,c}^\infty(M)/\mathcal W_{n-k+1,c}^\infty(M)\right)^*\\
& = \Gamma_c^\infty({\Val}_{n-k}^\infty(T M))^*.
\end{align*}

Passing to $G$-invariants and using Lemma \ref{lemma:dual_sections} yields 
\begin{displaymath}
\left[\mathcal W_k^{-\infty}(M)/\mathcal W_{k+1}^{-\infty}(M)\right]^G = \Gamma^\infty(\Val_k^{-\infty}(TM))^G.
\end{displaymath}

Denote by $\beta_k:\mathcal W_k^{-\infty}(M)^G\to \Gamma^\infty(\Val_k^{-\infty}(TM))^G$ the induced map.  

Fix some point $x_0 \in M$ and let $G_0$ be the stabilizer of $G$ at $x_0$. Since $M$ is isotropic, 
\begin{displaymath}
\beta_k\mu(x_0) \in \Val_k^{-\infty}(T_{x_0} M)^{G_0}=\Val_k^{-\infty}(\R^{p,q})^{\OO(p,q)}. 
\end{displaymath}

By \cite{bernig_faifman_opq}, we may find a linear combination $\tau:=a\mu_k+b\bar \mu_k$ of intrinsic volumes such that $\beta_k\tau(x_0)=\beta_k\mu(x_0)$. By $G$-invariance, it follows that $\beta_k\tau=\beta_k\mu$, i.e. $\mu-\tau \in \mathcal W_{k+1}^{-\infty}(M)^G$. Induction over $k$ then yields some $\tilde \tau \in \mathcal{LK}(M)$ with $\mu-\tilde \tau \in \mathcal{W}_{n+1}^{-\infty}(M)^G=\{0\}$. 

The argument in the case of generalized curvature measures is similar. Take $\Phi \in \mathcal{C}^{-\infty}(M)^G$, take $k \in \{0,\ldots,n+1\}$ with $\Phi \in \mathcal{C}_k^{-\infty}(M)$. If $k=n+1$, then $\Phi=0$ and we are done. Suppose that $k<n+1$. Fix a point $x_0 \in  M$. By Proposition \ref{prop_invariant_implies_horizontally_smooth}, $\Phi$ is horizontally smooth. Then $\sigma_k \Phi(x_0)  \in \Curv_k^{-\infty}(\R^{p,q})^{\OO(p,q)}$ is well-defined by Proposition \ref{prop_extension_solanes_wannerer_map}. By Proposition \ref{prop_flat_case} we find a linear combination $\Psi:=a \Lambda_k+b \bar \Lambda_k \in \widetilde{\mathcal{LK}}(M)$ of Lipschitz-Killing measures such that $\sigma_k \Phi(x_0)=\sigma_k \Psi(x_0)$. By $G$-invariance, we have $\sigma_k \Phi=\sigma_k\Psi$, hence $\Phi-\Psi \in \mathcal{C}_{k+1}^{-\infty}(M)^G$ and by induction on $k$ we find some $\tilde \Psi \in \widetilde{\mathcal{LK}}(M)$ with $\Phi-\tilde \Psi \in \mathcal{C}_{n+1}^{-\infty}(M)^G=\{0\}$.
\endproof

\def\cprime{$'$}


 \end{document}